\documentclass{tCON2e}

\usepackage{subfigure}
\usepackage{apalike}

\usepackage{enumerate}

\theoremstyle{plain}
\newtheorem{theorem}{Theorem}

\newtheorem{lemma}{Lemma}
\newtheorem{proposition}{Proposition}

\theoremstyle{definition}
\newtheorem{definition}{Definition}
\newtheorem{example}{Example}

\theoremstyle{remark}
\newtheorem{remark}{Remark}

\newcommand{\pfrac}[2]{\frac{\partial#1}{\partial #2}}
\newcommand{\mspan}{\mathrm{span}\,}
\newcommand{\mrk}{\mathrm{rk}\,}

\newcommand{\mmod}{\,\mathrm{mod}\,}

\newcommand{\dleq}[3]{#1 \leq #2 \leq #3}

\def \a {\alpha}
\def \b {\beta}

\def \S {\Sigma}

\def \R {\mathbb{R}}

\def \t {\tilde}

\def \vp {\varphi}

\def \cB {\mathcal{B}}

\newcommand*{\refname}{References}
\begin{document}



\title{Multi-input control-affine systems static feedback equivalent to a triangular form and their flatness}

\author{Shunjie Li$^{\rm a}$$^{\rm c}$ Florentina Nicolau$^{\rm b}$ and Witold  Respondek$^{\rm b}$$^{\rm c}$$^{\ast}$\thanks{$^\ast$Corresponding author. Email: witold.respondek@insa-rouen.fr \vspace{6pt}} \\
$^{a}${\em{School of mathematics and statistics, Nanjing University of Information Science and Technology, 210044,Nanjing, China}}\\
$^{b}${\em{Normandie Universit\'e, INSA de Rouen, Laboratoire de Math\'{e}matiques,76801, Saint-Etienne-du-Rouvray, France}}\\
\received{v4.0 released February 2014}
$^{c}${\em{Institute of Cyber-Systems and Control, Zhejiang University, 310027, Hangzhou, China}}
}

\maketitle

\begin{abstract}
In this paper, we give a complete geometric characterization of systems locally static feedback equivalent to a triangular form compatible with the chained form, for $m=1$, respectively with the $m$-chained form, for $m\geq2$. They are $x$-flat systems.  We provide a system of first order PDE's to be solved in order to find all $x$-flat outputs, for $m=1$, respectively all minimal  $x$-flat outputs, for $m\geq 2$. We illustrate our results by examples, in particular {by an application} to a mechanical system: the coin rolling without slipping on a moving table.
\end{abstract}

\section{Introduction}
The notion of flatness has been introduced in control theory in the 1990's by \cite{fliess1992systemes,fliess61vine},  see also \cite{isidori1986sufficient,jakubczyk1993invariants,martin1992phd,pomet1995differential}, and has attracted a lot of attention because of its multiple applications in the problem of trajectory tracking, motion planning and constructive controllability (see, e.g. \cite{fliess1999lie,levine2009analysis,martin2003flat,da2001flatness,pomet1997dynamic,respondek2003symmetries,schlacher2007construction}).

The fundamental property of flat systems is that all their solutions may be parameterized by $m$ functions and their time-derivatives, $m$ being the number of controls. More precisely, consider a nonlinear control  system
$$
\Xi\,:~\ \dot x=~F(x,u)
$$
where $x$ is the state defined on an open subset $X$ of  $\R^n$,  $u$ is the control taking values in an open subset $U$ of $\R^{m}$ (more generally, an $n$-dimensional manifold $ X$ and an $m$-dimensional manifold $U$, respectively) and the dynamics $F$ are smooth (the word smooth will always mean $\mathcal{C}^\infty$-smooth).
The system $\Xi$ is \textit{flat} if we can find $m$ functions, $\varphi_i(x,u,\dots,u^{(r)})$, for some $r \geq 0$, called \textit{flat outputs}, such that
\begin{equation}
\label{chap 3.1 descript}
x=\gamma(\varphi,\dots,\varphi^{(s)}) \mbox{ and }
u=\delta(\varphi,\dots,\varphi^{(s)}),
\end{equation}
for a certain integer $s$ and suitable maps $\gamma$ and $\delta$, where $\varphi=(\varphi_1,\dots,\varphi_m)$. Therefore all state  and control variables can be determined from the flat outputs without integration and all trajectories of the system can be completely parameterized. In the particular case $\varphi_i =\varphi_i(x)$, for $\dleq 1 i m$, we will say that the system is $x$-flat.
The minimal number of derivatives of components of a flat output $\varphi$, needed to express $x$ and $u$, will be called the differential weight of $\varphi$ {(see Section~\ref{chap 3.1 sec:flatness} for precise definitions).}
\bigskip

The problem of flatness of driftless two-input control-linear systems of the form
$$
\S_{lin}: \dot x = u_0g_0(x) + u_1g_1(x),
$$
defined on a open subset $X$ of {$\R^{n}$,} has been solved by \cite{martin1994feedback} (see also \cite{li2010flat,martin1993feedback} and a related result of \cite{cartan1914equivalence}). According to their result,  on an open and dense subset $X'$ of $X$, the system  $\S_{lin}$ is flat if and only if, its associated distribution $\mathcal{G} = \mspan \{g_0, g_1\}$ can be locally brought into the Goursat normal form, or equivalently, the  control system $\S_{lin}$ is locally static feedback equivalent to the chained form:

$$
\label{chap 3.1 Chained form}
Ch_1^k:\left\{
 \begin{array}{l lcl}
  \dot z_0 =v_0 & \dot z_1 &=&  z_2v_0 \\
                &  \dot z_2 &=&  z_3v_0 \\
                &     &\vdots&   \\
                &  \dot z_{k-1} &=& z_kv_0\\
                &   \dot z_{k} &=& v_1
 \end{array}\right.
$$
where $n=k+1$.

The first who noticed the existence  of singular points in the problem of transforming a distribution of rank two into the Goursat normal form  were \cite{giaro1978lecture}. \cite{murray1994nilpotent} presented a regularity condition that guarantees the feedback equivalence of $\S_{lin}$ to the chained form $Ch_1^k$ around an arbitrary point $x^*$. \cite{li2010flat} studied and solved the following problem: can a driftless two-input system be locally flat at a singular point of $\mathcal{G}$? In other words, can $\S_{lin}$ be flat without being locally equivalent to the chained form? Their result shows that a Goursat structure is $x$-flat only at regular points of $\mathcal{G}$. They also described all $x$-flat outputs and showed that they are parametrized by an arbitrary function of three variables canonically defined up to a diffemorphism.
\bigskip

In this paper we give a generalization of these results. Our goal is to characterize control-affine systems that are static feedback equivalent to the following triangular form
$$
 TCh_1^k:\left\{
 \begin{array}{l lcl c l}
  \dot z_0 =v_0 & \dot z_1 &=&  f_1(z_0,z_1,z_2) &+& z_2v_0 \\
                & \dot z_2 &=&  f_2(z_0,z_1,z_2, z_3)&+& z_3v_0 \\
                &  &\vdots&  &&  \\
                &  \dot z_{k-1} &=& f_2(z_0,\cdots, z_k)&+&z_kv_0\\
                &   \dot z_{k} &=& v_1 & &
 \end{array}\right.
$$
compatible with the chained form. Indeed, notice that  {in the $z$-coordinates} the distribution spanned by the controlled vector fields is in the chained form (Goursat normal form) and the drift has a triangular structure.

We will completely characterize control-affine systems that are static feedback equivalent to $TCh_1^k$ and show how their geometry differs and how it reminds that of  control-linear systems feedback equivalent to the chained form. Then, we will extend this result to the triangular form compatible with the $m$-chained form, i.e., we will characterize control-affine systems with $m+1$ inputs, where $m\geq 2$, that are static feedback equivalent to a normal form obtained by replacing {$z_j$, in $TCh_1^k$, by the vector $z^j = (z^j_1, \cdots, z^j_m)$, the smooth functions $f_j$ by $f^j = (f^j_1, \cdots, f^j_m)$} and the control $v_1$ by  the control vector $(v_1, \cdots, v_m)$. 
This form will be denoted by  $TCh_m^k$. Its associated distribution $\mathcal{G}= \mspan \{g_0, \cdots, g_m\}$, where $g_i$, for $\dleq 0 i m$, are the controlled vector fields, is called a Cartan distribution (or a contact distribution) for curves, see \cite{bryant1991exterior,olver1995equivalence,krasilchchik1986geometry}. 
The problem of characterizing {control-linear} systems that are locally static feedback equivalent to the $m$-chained form (or equivalently, that of characterizing Cartan distributions for curves) has been studied and solved (\cite{pasillas2001canonical}, see also \cite{mormul2002multi,pasillas2000geometry,pasillas2001contact,shibuya2009drapeau,yamaguchi1982contact}). It is immediate that systems locally feedback equivalent to the $m$-chained form are flat and in \cite{respondek2003symmetries}, all their minimal flat outputs (i.e., those whose differential weight is the lowest among all flat outputs of the system) have been described.

It is easy to see that the normal form $TCh_1^k$ (respectively $TCh_m^k$) is $x$-flat at any point of  $ X\times \mathbb{R}^{2}$ (respectively $ X\times \mathbb{R}^{m+1}$) satisfying some regularity conditions and we describe all its  $x$-flat outputs (respectively all its minimal $x$-flat outputs). Their description reminds very much that of control-linear systems feedback equivalent to the chained form, for $m=1$, respectively to the $m$-chained form, for $m\geq 2$, although new phenomena appear related to singularities in the state and control-space. 

{Since $ TCh_1^k$ and$TCh_m^k$ are flat, the} paper gives sufficient conditions for a system to be $x$-flat. We will also show that these conditions are not necessary for $x$-flatness of control-affine system whose associated distribution spanned by the controlled vector fields $\mathcal{G}= \mspan \{g_0, \cdots, g_m\}$ is feedback equivalent to the $m$-chained form. Indeed, we show that there are $x$-flat control-affine systems for which there exist local coordinates in which the distribution spanned by the controlled vector fields has the $m$-chained structure but the drift is not triangular (see Example~\ref{chap 3.1 ssec_academical example}).

The triangular form $TCh^k_1$ was considered in \cite{li2013characterization}, where its flatness was observed but its description was not addressed. A characterization of $TCh^k_1$ has been recently proven by  \cite{silveira2010formas} and by  \cite{silveira2013flat}, where a solution dual to ours (using an approach based on differential forms and codistributions rather than distributions) is given.
Our aim is to treat in a homogeneous way the two-input case of $TCh^k_1$ and the multi-input case of
$TCh^k_m$, using the formalism of vector fields and distributions, as well as to describe all flat outputs and their singularities (which are more natural to deal with in the language of vector fields).

\bigskip

The paper is organized as follows. In Section \ref{chap 3.1 sec:flatness}, we recall the definition of flatness and define the notion of  differential weight of a flat system. In Section \ref{chap 3.1 sec:main-results}, we give our main results: we characterize  control-affine systems static feedback equivalent to the  triangular form $TCh_1^k$, for $m=1$, and to $TCh_m^k$, for $m\geq 2$. 
We describe in Section~\ref{chap 3.1 sec:Flat-outputs} all minimal flat outputs including their singularities and we study  also singular control values at which the system ceases to be flat. 
Moreover, we give also in that  section a system of first order PDE's to be solved in order to find all  $x$-flat outputs, for $m=1$, and all minimal $x$-flat outputs, for $m\geq 2$. We illustrate our results by two examples in Section~\ref{chap 3.1 sec:application} and provide proofs in Section~\ref{chap 3.1 sec:proofs}.

\section{Flatness} \label{chap 3.1 sec:flatness}

Fix an integer $l \geq -1$ and denote $ U^l = U \times \mathbb{R}^{ml}$ and $\bar u^l = (u,\dot  u, \dots, u^{(l)})$. For $l=-1$, the set $ U^{-1}$ is empty and $\bar u^{-1}$ is an empty sequence.

\smallskip

\begin{definition} \label{chap 3.1 Def:flatness}
The system $\Xi: \dot x=~F(x,u)$ is \textit{flat} at $(x^*, \bar {u}^{*l}) \in X \times U^{l}$, for  $l\geq -1$, if there exists a neighborhood $\mathcal O^{l}$ of $(x^*, \bar {u}^{*l})$
and $m$ smooth functions $\varphi_i=\varphi_i(x, u,\dot  u, \dots, u^{(l)})$, $1 \leq i\leq m$, defined in $\mathcal O^l$, having the following property: there exist an integer $s$ and smooth functions $ \gamma_i$, $1 \leq i \leq  n$, and $\delta_j$, $ 1 \leq j \leq m$, such that
$$x_i = \gamma_i (\varphi, \dot \varphi, \dots, \varphi^{(s)}) \mbox{ and }u_j=  \delta_j (\varphi, \dot \varphi, \dots, \varphi^{(s)})
$$
along any  trajectory $x(t)$ given by a control $u(t)$ that satisfy $(x(t), u(t), \dots , u^{(l)}(t)) \in \mathcal O^l$, where $\varphi=(\varphi_1, \dots, \varphi_m)$ and is  called \textit{flat output}.
\end{definition}

When necessary to indicate the number of derivatives of $u$ on which the flat outputs $\varphi_i$ depend, we will say that the system $\Xi$ is $(x, u, \cdots, u^{(r)})$-flat if $ u^{(r)}$ is the highest derivative on which $\varphi_i$ depend and in the particular case $\varphi_i =\varphi_i(x)$, we will say that the system is $x$-flat.
In general, $r$ is smaller than the integer $l$ needed to define the neighborhood $\mathcal O^l$ which, in turn, is smaller than the number of derivatives of $\varphi_i$ that are involved. In our study, $r$ is always equal to~-1, i.e., the flat outputs depend on $x$ only, and  $l$ is 0.

The minimal number of derivatives of components of a flat output $\varphi$, needed to express $x$ and~$u$, will be called the  differential weight of that flat output and will be formalized as follows.
By definition, for  any flat output $\varphi$ of $\Xi$ there exist integers $s_1,\dots,s_m$ such that
$$
\left.
\begin{array}{lll}
x&= &\gamma(\varphi_1,\dot\varphi_1,\dots,\varphi_1^{(s_1)},\dots, \varphi_m,\dot\varphi_m, \dots, \varphi_m^{(s_m)}) \vspace{0.2cm}\\
u&= & \delta(\varphi_1,\dot\varphi_1,\dots,\varphi_1^{(s_1)},\dots, \varphi_m,\dot\varphi_m, \dots, \varphi_m^{(s_m)}),
\end{array}
\right.
$$
Moreover, we can choose $(s_1,\dots,s_m)$ such that (see \cite{respondek2003symmetries}) if for any other {$m$}-tuple $(\t s_1,\dots,\t s_m)$ we have
$$
\left.
\begin{array}{lll}
x&= &\t \gamma(\varphi_1,\dot\varphi_1,\dots,\varphi_1^{(\t s_1)},\dots, \varphi_m,\dot\varphi_m, \dots, \varphi_m^{(\t s_m)})  \vspace{0.2cm}\\
u&= &\t \delta(\varphi_1,\dot\varphi_1,\dots,\varphi_1^{(\t s_1)},\dots, \varphi_m,\dot\varphi_m, \dots, \varphi_m^{(\t s_m)}),
\end{array}
\right.
$$
then $s_i \leq \t s_i$, for $\dleq 1 i m$.

We will call  $\sum_{i=1}^{m}(s_i+1) =m+\sum_{i=1}^{m}s_i $ the differential weight of $\varphi$.  A flat output of $\Xi$ is called \textit{minimal} if its differential weight is the lowest among all flat outputs of $\Xi$. We define the \textit{differential weight} of a flat system to be equal to the differential weight of a minimal flat output.
\bigskip

\section{Main results: characterization of the triangular form}\label{chap 3.1 sec:main-results}

From now on, we will denote the number of controls by $m+1$ (and not by $m$) since, as we will see below, for all classes of systems that follow one control plays a particular role.

Consider {the} control-affine system
\begin{equation}\label{chap 3.1 Sigma}
 \Sigma_{aff}: \dot x = f(x) + \sum_{i=0}^m u_ig_i(x),
\end{equation}
defined on  an open subset $X$ of $\mathbb{R}^{ n}$, where $n=km+1$ (or an $n$-dimensional manifold $X$), where~$f$ and $g_0, \cdots, g_m$ are smooth vector fields on $X$ and the number of controls is $m+1\geq 2$.

To $\Sigma_{aff}$ we associate the following distribution $\mathcal{G} = \mspan \{g_0, \cdots, g_m\}$. We define inductively the derived flag  of $\mathcal{G}$ by
$$
\mathcal{G}^0 = \mathcal{G} \mbox{ and } \mathcal{G}^{i+1} =  \mathcal{G}^i + [ \mathcal{G}^i,  \mathcal{G}^i],\, i \geq 0.
$$

Let $\mathcal{D}$ be a non involutive distribution of rank $d$, defined on $X$ and define its annihilator $\mathcal{D}^{\perp}=\{\omega\in \Lambda^1(X)\,:\,<\omega,f >=0, \forall f\in\mathcal{D}\}$, where $ \Lambda^1(X)$ stands for the collection of smooth differential 1-forms on $X$. A vector field $c\in\mathcal{D}$ is called characteristic for $\mathcal{D}$ if it satisfies $[c,\mathcal{D}] \subset \mathcal{D}$. The characteristic distribution of $\mathcal{D}$, denoted by $\mathcal{C}$, is the distribution spanned by all its characteristic vector fields, i.e., 
$$\mathcal{C} = \{c \in \mathcal{D}: [c,\mathcal{D}]\subset \mathcal{D}\}$$
and can be computed as follows.
Let $\omega_1, \dots, \omega_q$, where $q= n-d$, be differential 1-forms locally spanning the annihilator of $\mathcal{D}$, that is $\mathcal{D}^{\perp}= \mspan \{\omega_1, \dots, \omega_q\}$.  For any $\omega \in \mathcal{D}^{\perp}$, we define $\mathcal{W}(\omega) = \{f \in \mathcal{D}: f \lrcorner\, d\omega\in \mathcal{D}^\perp\}$, where $\lrcorner$ is the interior product. The characteristic distribution of $\mathcal{D}$ is given by
$$
\mathcal{C} = {\bigcap}_{i=1}^q \mathcal{W}(\omega_i).
$$
It follows directly from the Jacobi identity that the characteristic distribution is always involutive.
 \bigskip

Our main results describing control-affine systems locally static feedback equivalent to the triangular form compatible to the chained form and to the $m$-chained form, are given by the two following theorems corresponding to two-input control-affine systems, i.e., $m=1$ (Theorem \ref{chap 3.1 thm:m=1}), and to control-affine systems with $m+1$ inputs, for $m\geq 2$ (Theorem \ref{chap 3.1 thm:m>=2}). Let us first consider the case $m=1$, which has also been solved, using the formalism of differential forms and codistributions, by  \cite{silveira2010formas} and by \cite{silveira2013flat}.

\begin{theorem}\label{chap 3.1 thm:m=1}
Consider a two-input control-affine system $\Sigma_{aff}$, given by (\ref{chap 3.1 Sigma}), for $m=1$, and fix $x^* \in X$, an open subset of $\mathbb{R}^{k+1}$.
 The system $\Sigma$ is locally, around $x^*$, static feedback equivalent to the triangular form $TCh_1^k$ if and only if the following conditions are satisfied:
    \begin{enumerate}[{(Ch}1{)}]
    \item $\mathcal{G}^{k-1} = TX$;
    \item $\mathcal{G}^{k-3}$ is of constant rank $k-1$ and, moreover,  the characteristic distribution $\mathcal{C}^{k-2}$ of $\mathcal{G}^{k-2}$ is contained in $\mathcal{G}^{k-3}$ and has constant corank one in $\mathcal{G}^{k-3}$;
    \item $\mathcal{G}^0(x^*)$ is not contained in $\mathcal{C}^{k-2}(x^*)$;
    \end{enumerate}

    \begin{enumerate}[(Comp)]
    \item $[f, \mathcal{C}^i] \subset  \mathcal{G}^i$, for $\dleq 1 i k-2$, where $\mathcal{C}^i$ is the characteristic distribution of $\mathcal{G}^i$.
    \end{enumerate}
\end{theorem}

It was stated and proved in \cite{pasillas2001canonical} that items \textit{(Ch1)-(Ch3)} characterize, locally, the chained form (or equivalently  the Goursat normal form). Therefore, they are equivalent to the well known conditions describing the chained form \cite{murray1994nilpotent} (see also \cite{kumpera1982equivalence,martin1994feedback,montgomery2001geometric,mormul2000goursat,pasillas2001geometry}):
\textit{\begin{enumerate}[{(Ch}1{)'}]
 \item $\mrk \mathcal{G}^i= i+2$, for $\dleq 0 i k-1$,
 \item $\mrk \mathcal{G}^i(x^*) = \mrk \mathcal{G}_i(x^*) = i+2$, for $\dleq 0 i k-1$,
 where the distributions $ \mathcal{G}_i$ form the Lie flag of $\mathcal{G}$ and are defined by
$\mathcal{G}_0 = \mathcal{G} \mbox{ and } \mathcal{G}_{i+1} =  \mathcal{G}_i + [ \mathcal{G}_0,  \mathcal{G}_i],\, i \geq 0,$
\end{enumerate}}
\noindent
and assure the existence of a change of coordinates $z=\phi(x)$ and of an invertible static feedback transformation of the form $u=\beta \t u$, bringing the control vector fields $g_{0}$ and $g_{1}$ into the chained form. 

Item $(Comp)$ takes into account the drift and gives the compatibility conditions for $f$ to have the desired triangular form in the right system of coordinates, i.e., in coordinates $z$ in which the controlled vector fields are in the chained form.

Since the distribution $\mathcal{G}$, associated to $\Sigma_{aff}$, satisfies $(Ch1)'$, all characteristic distributions  $\mathcal{C}^i$ of  $\mathcal{G}^i$ are well defined, for $\dleq 1 i k-2$. Indeed, recall the following result due to  \cite{cartan1914equivalence}:

 \begin{lemma} \emph{(E. Cartan)} \label{chap 3.1 lemma Cartan}
 Consider a rank two distribution $\mathcal{G}$ defined on a manifold $X$ of dimension $k + 1$, for $k \geq 3$. If $\mathcal{G}$ satisfies $\mrk \mathcal{G}^i= i+2$, for $\dleq 0 i k-1$, everywhere on~$X$, then each distribution $\mathcal{G}^i$, for $\dleq 0 i k-3$, contains a unique involutive subdistribution $\mathcal{C}^{i+1}$ that is characteristic for $\mathcal G^{i+1}$ and has constant corank one in $\mathcal{G}^i$.
\end{lemma}

The conditions of the above theorem are verifiable, i.e., given a two-input control-affine system and an initial point $x^*$, we can verify whether it is locally static feedback equivalent, around $x^*$, to $TCh_1^k$ and verification (in terms of vector fields of the initial system) involves derivations and algebraic operations only, without solving PDE's.
\bigskip

Next, we consider the case $m\geq 2$ and extend the above result to a triangular form compatible with the $m$-chained form. An $(m+1)$-input driftless control system $\S_{lin}: \dot z = \sum_{i=0}^mv_ig_i(z)$, defined on $\mathbb{R}^{km+1}$, is said to be in the $m$-chained form if it is represented by

$$
Ch_m^k:\left\{
\begin{array}{l lcl c lcl}
 \dot z_0 =v_0 & \dot z_1^1&=&  z_1^2v_0& \cdots & \dot z_m^1 &=&  z_m^2v_0\vspace{0.2cm}\\
                 & \dot z_1^2 &=& z_1^3v_0&   & \dot z_m^2 &=&   z_m^3v_0\\
                 &            &\vdots&    & & & \vdots& \\
                 & \dot z_1^{k-1} &=&  z_1^kv_0&  { \cdots} & \dot z_m^{k-1} &=&   z_m^kv_0\vspace{0.2cm}\\
                 & \dot z_1^{k} &=& v_1& { \cdots} & \dot z_m^{k} &=& v_m

\end{array}
\right.
$$
\smallskip

Denote {$\bar z^j = (z_1^1, \cdots z_m^1, z_1^2, \cdots z_m^2,\cdots,  z_1^j, \cdots z_m^j)$, for $\dleq 2 j k$}. Our goal is to  characterize the following triangular normal form
$$
TCh_{m}^k:\left\{
\begin{array}{l lclcl c lclcl}
 \dot z_0 =v_0 & \dot z_1^1&=& f_1^1( z_0,\bar z^2) & + & z_1^2v_0& \cdots & \dot z_m^1 &=& f_m^1( z_0,\bar z^2)& + &  z_m^2v_0\vspace{0.2cm}\\
                 & \dot z_1^2 &=& f_1^2( z_0,\bar z^3)  & + &  z_1^3v_0&   & \dot z_m^2 &=& f_m^2( z_0,\bar z^3) & + &  z_m^3v_0\\
                 &            &\vdots&  & &  & & & \vdots& \\
                 & \dot z_1^{k-1} &=& f_1^{k-1}( z_0,\bar z^k)  & + &  z_1^kv_0&  \cdots & \dot z_m^{k-1} &=& f_m^{k-1}( z_0,\bar z^k) & + &  z_m^kv_0\vspace{0.2cm}\\
                 & \dot z_1^{k} &=& v_1& &&  \cdots & \dot z_m^{k} &=& v_m

\end{array}
\right.
$$
with $m+1$ inputs, $m\geq 2$.
Theorem \ref{chap 3.1 thm:m>=2} below gives necessary and sufficient conditions for a control system to be locally static feedback equivalent to  $TCh_m^k$.

\begin{theorem}\label{chap 3.1 thm:m>=2}
Consider a control-affine system $\Sigma_{aff}$, given by (\ref{chap 3.1 Sigma}),, on an open subset $X$ of $\mathbb{R}^{km+1}$, for $m\geq 2$, and fix $x^* \in X$. The system $\Sigma_{aff}$ is locally, around $x^*$, static feedback equivalent to the triangular form $TCh_{m}^k$ if and only if the following conditions are satisfied:
    \begin{enumerate}[{(m-Ch}1{)}]
    \item $\mathcal{G}^{k-1} = TX$;
    \item $\mathcal{G}^{k-2}$ is of constant rank $(k-2)m+1$ and contains an involutive subdistribution $\mathcal{L} $ that has constant corank one in $\mathcal{G}^{k-2}$;
    \item $\mathcal{G}^0(x^*)$ is not contained in $\mathcal{L}(x^*)$;
    \end{enumerate}
    \begin{enumerate}[(m-Comp)]
    \item $[f, \mathcal{C}^i] \subset  \mathcal{G}^i$, for $\dleq 1 i k-2$, where $\mathcal{C}^i$ is the characteristic distribution of $\mathcal{G}^i$.
    \end{enumerate}
\end{theorem}

In order to verify the conditions of Theorem \ref{chap 3.1 thm:m>=2}, we have to check whether the distribution $\mathcal{G}^{k-2}$ contains an involutive subdistribution $\mathcal{L}$ of corank one. Checkable necessary and sufficient conditions for the existence of such an involutive subdistribution, together with a construction, follow from the work of  \cite{bryant1979some} and are given explicitly in \cite{pasillas2001contact}.  We present in Appendix A 
the conditions for the existence and construction {of $\mathcal{L}$}. In our case, if such a distribution exists, it is always unique. As a consequence, all conditions of Theorem \ref{chap 3.1 thm:m>=2} are verifiable, i.e., given a control-affine system and an initial point $x^*$, we can verify whether it is locally static feedback equivalent, around $x^*$, to $TCh_m^k$ and verification involves derivations and algebraic operations only, without solving PDE's.   

Conditions \textit{(m-Ch1)-(m-Ch3)} characterize the $m$-chained form \cite{pasillas2001canonical} (see also \cite{pasillas2000geometry,pasillas2001contact}) and assure the existence of a change of coordinates $z=\phi(x)$ and of an invertible static feedback transformation of the form $u=\beta \t u$, bringing the control vector fields $g_{i}$ into the $m$-chained form. We define the diffeomorphism $\phi$ {and the feedback transformation $\beta$} in Appendix B. {The diffemorphism  $\phi$  defines also the coordinates in which the system takes  the triangular form $TCh_m^k$.}

Item \textit{(m-Comp)} takes into account the drift and gives the compatibility conditions for $f$ to have the desired triangular form in the right system of coordinates, i.e., in $z$-coordinates  in which the controlled vector fields are in the $m$-chained form. {Formally it has the same form as $(Comp)$ in the case $m=1$.}

The characteristic distributions  $\mathcal{C}^i$, for $\dleq 1 i k-2$, are well defined and have corank one in~$\mathcal{G}^{i-1}$. Indeed, recall the following result stated in \cite{pasillas2001canonical}:

\begin{lemma}
{Assume that a distribution $\mathcal{G}$ defined on a manifold $X$ of dimension $km + 1$ satisfies the conditions \textit{(m-Ch1)-(m-Ch3)} of Theorem~\ref{chap 3.1 thm:m>=2}. Then} $\mathcal{G}^i$ has constant rank $(i+1)m+1$, for $\dleq 0 i k-2$, and contains an involutive subdistribution $\mathcal{L}^i$ of corank one in $\mathcal{G}^i$. Moreover  $\mathcal{L}^i$ is the unique corank one subdistribution satisfying this property, for $\dleq 0 i k-2$, and it coincides with the characteristic distribution $\mathcal{C}^{i+1}$ of  $\mathcal G^{i+1}$, for $\dleq 0 i k-3$.
%
\end{lemma}

It has been shown in \cite{respondek2001transforming} (see also \cite{pasillas2001canonical}) that all information about the distribution~$\mathcal{G}$ is encoded completely in the existence of the last involutive subdistribution $\mathcal{L}^{k-2}$ (being, actually, the involutive distribution $\mathcal{L}$ of item \textit{(m-Ch2)} of Theorem~\ref{chap 3.1 thm:m>=2}) which implies the existence of all involutive subdistributions $\mathcal{L}^i =\mathcal{C}^{i+1}$, for $\dleq 0 i k-3$.
\bigskip

The characterization of the chained form (conditions \textit{(Ch1)-(Ch3)} of  Theorem \ref{chap 3.1 thm:m=1}) and that of the $m$-chained form (\textit{(m-Ch1)-(C-mCh3)} of Theorem \ref{chap 3.1 thm:m>=2}) are different, but compatibility conditions are the same, compare \textit{(Comp)} and \textit{(m-Comp)}. The involutive subdistribution $\mathcal{L}$, which is crucial for the $m$-chained form, is absent in the compatibility conditions, but plays a very important role in calculating minimal flat outputs and in describing singularities (see Section~ \ref{chap 3.1 sec:Flat-outputs}).

\section{Flatness and flat outputs description}\label{chap 3.1 sec:Flat-outputs}

In this section, firstly, we discuss flatness of control systems static feedback equivalent to $TCh_1^k$, respectively to $TCh_m^k$. Secondly, we answer the question whether a given pair   (respectively an $(m+1)$-tuple) of smooth functions on $X$ is an $x$-flat output for a system static feedback equivalent to $TCh_1^k$ (respectively a minimal $x$-flat output for a system static feedback equivalent to $TCh_m^k$) and, finally, provide a system of PDS's to be solved in order to find all these flat outputs. In particular, we will discuss their uniqueness, their singularities, and compare their description with that of flat outputs for the chained form (respectively for the $m$-chained form).

{\subsection{Flatness of control systems static feedback equivalent to $TCh_1^k$}} \label{chap 3.1 ssec_flatness m=1}

Let us first consider the case $m=1$.
It is clear that $TCh_1^k$ is $x$-flat, with $\varphi = (z_0,z_1)$ being a flat output around any point  $(z^*,v^*)$ satisfying
$$
\pfrac{f_i}{z_{i+1}}(z^*)+v^*_{0}\neq 0,  \mbox{ for } 1 \leq i \leq k-1,
$$
where $v^*=(v^*_{0},v^*_{1}).$
Therefore control systems equivalent to $TCh_1^k$  are $x$-flat and exhibit a singularity in the control space (depending on the state) which we will describe in an invariant way as follows. For $\mathcal{C}^{1}\subset \mathcal{C}^{2}\subset \cdots \subset \mathcal{C}^{k-2}$, the sequence of characteristic distributions $\mathcal{C}^{i}$ of $\mathcal{G}^{i}$, for $1\leq i \leq k-2$, see Lemma \ref{chap 3.1 lemma Cartan}, choose vector fields $c_{1},\ldots,c_{k-2}$ such that $\mathcal{C}^{i}= \textrm{span}\,\{c_{1},\ldots,c_{i}\}$. For each $0\leq i \leq k-3$, define
$$
U_{sing}^{i}(x)=\left\{u^{i}(x)=(u^{i}_{0}(x),u^{i}_{1}(x))^{\top}:[f+u^{i}_{0}g_{0}+u^{i}_{1}g_{1}, \mathcal{C}^{i+1}]\subset \mathcal{G}^{i}\right\}.
$$
The controls $u^{i}(x)$ exist, are smooth, and for any $0\leq i \leq k-3$ define (for any fixed $x\in X$) a 1-dimensional affine subspace of $U=\mathbb{R}^{2}$. To see those three properties, notice that $[f,c_{i+1}]$, $[g_{0},c_{i+1}]$, and $[g_{1},c_{i+1}]$ span a distribution of rank one modulo~$\mathcal{G}^{i}$ (since all three belong to $\mathcal{G}^{i+1}$ and $\textrm{corank}(\mathcal{G}^{i}\subset\mathcal{G}^{i+1})= 1$) and either $[g_{0},c_{i+1}]$ or $[g_{1},c_{i+1}]$ (or both) does not vanish modulo~$\mathcal{G}^{i}$. 
To calculate $U_{sing}^{i}(x)$ explicitly, assume that we have chosen $(g_{0},g_{1})$ such that $g_{1}=c_{1}$. 
Then  {$[g_1,c_{i+1}]=$} $[c_1,c_{i+1}]\in \mathcal{G}^{i} $ and $[f,c_{i+1}]=\alpha[g_{0},c_{i+1}] \, \textrm{mod}\, \mathcal{G}^{i}$, for some smooth function $\alpha$. We put $u_{0}^{i}(x)=-\alpha(x)$ and $u_{1}^{i}(x)$ arbitrary. It is clear that the definition of $(u_{0}^{i}(x),u_{1}^{i}(x))$ does not depend on the choice of $c_{1},\ldots,c_{k-2}$ and is feedback invariant (independently of whether we have chosen $g_{1}=c_{1}$ or not).
Indeed, if $u^{i}(x)\in U_{sing}^{i}(x)$, then for the feedback modified system $\dot{x}=\tilde{f}+\tilde{g}\tilde{u}$, where $\tilde{f}=f+g\alpha$ and $\tilde{g}=g\beta$, it is the feedback modified control $\tilde{u}^{i}=\beta^{-1}(-\alpha+u^{i})$ that, clearly, satisfies $\tilde{u}^{i}\in U_{sing}^{i}$.

Let $\mathcal{L}$ be any involutive distribution of corank two in $TX$ such that $\mathcal{L}\subset \mathcal{G}^{k-2}$. 
Fix $l\in \mathcal{L}$ such that $l \not\in \mathcal{C}^{k-2}$ and put
$$
U_{_{\mathcal{L}-sing}}^{k-2}(x)= \left\{u^{k-2}(x) = (u_{0}^{k-2}(x),u_{1}^{k-2}(x))^{\top}:[f+u_{0}^{k-2}g_{0}+u_{1}^{k-2}g_{1},l]\in \mathcal{G}^{k-2}\right\}.
$$
If $\mathcal{G}^{0}(x^{\ast}) \not\subset \mathcal{L}(x^{\ast})$, where $x^{\ast}$ is a nominal point around which we work, then  the controls $u^{k-2}(x)$ exist, are smooth, and (for any fixed $x\in X$) form a 1-dimensional affine subset of $U=\mathbb{R}^{2}$ because $\mathcal{G}^{k-2}$ is of corank one in $TX$ and either $[g_{0},l]$ or $[g_{1},l]$ is not in $\mathcal{G}^{k-2}$. 
If $\mathcal{G}^{0}(x^{\ast}) \subset \mathcal{L}(x^{\ast})$, then under the assumption, which we will always assume, $( d \varphi_{0} \wedge  d \varphi_{1} \wedge  d \dot{\varphi}_{0} \wedge d \dot{\varphi}_{1})(x^{\ast}, u^{\ast}) \neq 0$, where the functions $\varphi_{0}$ and $\varphi_{1} $ are such that $\mathcal{L}^\perp = \mspan\{{d}\varphi_{0},{d}\varphi_{1} \}$, we have $u^{\ast}\not\in U_{_{\mathcal{L}-sing}}^{k-2}(x^*) $ and in $\mathcal{X}^*\times \mathbb{R}^2$, where $\mathcal{X}^* $ is a sufficiently small neighborhood of $x^*$, the set $U_{_{\mathcal{L}-sing}}^{k-2}(x) $ consists of two connected components that define, for each fixed value $x \in \mathcal{X}^*$, $x\neq x^*$, an affine subspace of $U= \mathbb{R}^2$.

Clearly $U_{_{\mathcal{L}-sing}}^{k-2}$ is feedback invariant and does not depend on the choice of $l\in \mathcal{L}$ but it depends on the distribution $\mathcal{L}$. Define
$$
U_{sing}^{k-2}=\bigcap_{\mathcal{L}} U_{_{\mathcal{L}-sing}}^{k-2}
$$
where the intersection is taken over all $\mathcal{L}$ as above, that is,  involutive distribution of corank two in $TX$, satisfying $\mathcal{L}\subset \mathcal{G}^{k-2}$ . 
Define
$$
U_{sing}=\bigcup_{i=0}^{k-3}U_{sing}^{i} \cup U_{sing}^{k-2}
$$
and
$$
U_{_{\mathcal{L}-sing}}=\bigcup_{i=0}^{k-3}U_{sing}^{i} \cup U_{_{\mathcal{L}-sing}}^{k-2}.
$$
We will use both sets in Theorem \ref{chap 3.1 thm:output:m=1} describing controls singular for flatness and in Proposition~\ref{chap 3.1 prop output m=1 S_lin} comparing flat outputs of the triangular form $TCh_1^k$ with those of the associated chained form $Ch_1^k $.

\begin{theorem}\label{chap 3.1 thm:output:m=1}
Consider a two-input control-affine system $\Sigma_{aff}: \dot x = f(x) + u_0g_0(x) + u_1g_1(x)$, defined on an open subset $X$ of $\mathbb{R}^{k+1}$, where $k+1\geq 4$. 
Assume that $\Sigma_{aff}$ is locally, around $x^* \in X$, static feedback equivalent to $TCh_1^k$. Then we have:
\begin{enumerate} [(F1)]
\item $\Sigma_{aff}$ is $x$-flat at any $(x^*,u^*)\in X\times \mathbb{R}^2$ such that  $u^{*} \not\in U_{sing}(x^{*})$.
\item Let  $\vp_0$, $\vp_1$ be two smooth functions defined in a neighborhood $\mathcal{X}$ of $x^*$ and $g$ be an arbitrary vector field in $\mathcal{G}$ such that $g(x^*)\not\in \mathcal{C}^{k-2}(x^*)$. Then the following conditions are equivalent in~$\mathcal{X}$:

\begin{enumerate}[{(}i{)}]
\item The pair $(\vp_0,\vp_1)$ is an $x$-flat output of $\Sigma_{aff}$ at $(x^*,u^*) \in \mathcal{X}^*\times \mathbb{R}^2$, where  $\mathcal{X}^* $ is a neighborhood of $x^*$ ;

\item The pair $(\vp_0,\vp_1)$  satisfies the following conditions:
    \begin{enumerate}[{(FO}1{)}]
		\item $( d \varphi_{0} \wedge  d \varphi_{1} \wedge  d \dot{\varphi}_{0} \wedge d \dot{\varphi}_{1})(x^{\ast}, u^*) \neq 0$, where $\dot{\varphi}_{i}=L_{F_{aff}}\varphi_i$, for $i=0,1$ and $F_{aff}=f+u_0g_0+u_1g_1$;
		\item  $L_c \varphi_0 =  L_c\varphi_1=0$ and {$(L_g \varphi_0) (L_{[c,g]} \varphi_1) - (L_g \varphi_1)( L_{[c,g]} \varphi_0) =0$},
		for any $c \in \mathcal C^{k-2}$;
		 \item     {$u^* \not\in U_{_{\mathcal{L}-sing}}(x^*)$, where $\mathcal{L}=(\emph{span}\,\{ d \varphi_{0}, d \varphi_{1}\})^{\bot}$.}
    \end{enumerate}

 \item The pair $(\vp_0,\vp_1)$  satisfies the following conditions:
    \begin{enumerate}[{(FO}1{)$^{\prime}$}]
		\item $( d \varphi_{0} \wedge  d \varphi_{1} \wedge  d \dot{\varphi}_{0} \wedge d \dot{\varphi}_{1})(x^{\ast}, u^*) \neq 0$, where $\dot{\varphi}_{i}=L_{F_{aff}}\varphi_i$, for $i=0,1$, and $F_{aff}=f+u_0g_0+u_1g_1$;
		\item $\mathcal{L} = (\mspan \{ d \varphi_0, d \varphi_1\})^\perp\subset \mathcal{G}^{k-2}$;
		 \item $u^* \not\in U_{_{\mathcal{L}-sing}}(x^*)$. 
		\end{enumerate}

\end{enumerate}
\end{enumerate}

\end{theorem}
Notice that since $\S_{aff}$ is locally, around~$x^*$, static feedback equivalent to $TCh_1^k$, its associated control-linear system $\S_{lin}: \dot x=  u_0 g_0(x) + u_1g_1(x)$ is locally, around $x^*$, static feedback equivalent to the chained form $Ch_1^k$.
The {next} result shows how the similarities and differences between two-input control-linear systems and control-affine systems locally equivalent to $TCh_1^k$ are reflected by their flatness.
It turns out that flat outputs of $\S_{lin}$ are flat outputs of $\Sigma_{aff}$ (independently of the choice of $f$ although singular control values depend on $f$) and most of flat outputs of $\Sigma_{aff}$ are flat outputs of the corresponding $\S_{lin}$ but not all, as the following proposition explains. Define
$$
U_{char}(x)=\left\{u(x)=(u_{0}(x),u_{1}(x))^{\top}:(u_{0}g_{0}+u_{1}g_{1})(x)\in \mathcal{C}^{1}(x)\right\}.
$$

\begin{proposition}
\label{chap 3.1 prop output m=1 S_lin}
Consider a two-input control-affine system $\Sigma_{aff}: \dot x = f(x) + u_0g_0(x) + u_1g_1(x)$, defined on an open subset $X$ of $\mathbb{R}^{k+1}$, where $k+1\geq 4$, and its associated control-linear system $\S_{lin}: \dot x=  u_0 g_0(x) + u_1g_1(x)$.  Assume that $\Sigma_{aff}$ is locally, around $x^* \in X$, static feedback equivalent to $TCh_1^k$. Then we have:
\begin{enumerate} [(F1)]
\setcounter{enumi}{2}
\item $\Sigma_{lin}$ is $x$-flat at any $(x^*,u^*)\in X\times \mathbb{R}^2$ such that  $u^{*} \not\in U_{char}(x^{*})$.

\item A pair  $(\vp_0$, $\vp_1)$ of smooth functions defined in a neighborhood $\mathcal{X}$ of $x^*$
is an $x$-flat output of $\Sigma_{lin}$ at $(x^*,u^*) \in \mathcal{X^*}\times \mathbb{R}^2$ such that
{$\mathcal{X}^*\subset X$ is an open} neighborhood of $x^*$ and $u^* \not\in U_{char}(x^*)$ if and only if it satisfies the conditions (FO1)-(FO2) or, equivalently, (FO1)'-(FO2)' of Theorem \ref{chap 3.1 thm:output:m=1}, where  $\dot{\varphi}_{i}$, for $i=0,1$, is understood as
    $\dot{\varphi}_{i}=L_{F_{lin}\varphi_{i}}$ and $F_{lin}=u_0g_0+u_1g_1$;

\item If $(\vp_0,\vp_1)$ is a flat output of $\S_{lin}$ at $(x^*,u^*)$, where $u^* \not\in U_{char}(x^*)$, then $(\vp_0,\vp_1)$ is a flat output of $\S_{aff}$ at$(x^*,{\tilde u^*})$, where ${\tilde u^*} \not\in U_{_{\mathcal{L}-sing}}(x^*)$  with $\mathcal{L} = (\mspan\{d\vp_0, d\vp_1\})^\perp$.

\item Let $g$ be an arbitrary vector field in $\mathcal{G}$ such that $g(x^*)\not\in \mathcal{C}^{k-2}(x^*)$. If $(\vp_0,\vp_1)$ is a flat output of $\S_{aff}$ at $(x^*,{\tilde u^*})$, where ${\tilde u^*} \not\in U_{_{\mathcal{L}-sing}}(x^*)$,  with $\mathcal{L} = (\mspan\{d\vp_0, d\vp_1\})^\perp$, and satisfies $(L_g\vp_0,L_g\vp_1)(x^*)\neq (0,0)$,
    then $(\vp_0,\vp_1)$ is a flat output of $\S_{lin}$ at $(x^*,u^*)$,
    where $u^* \not\in U_{char}(x^*)$.

\end{enumerate}

\end{proposition}
For a pair of functions $(\vp_0,\vp_1)$, the conditions to be a flat output are, formally, the same for
$\S_{aff}$ and the associated control-linear system $\S_{lin}$ and are given by (FO1)-(FO2) (or, equivalently, by (FO1)'-(FO2)'). Notice, however, that the vector field along which we differentiate changes from
$F_{aff}$ into $F_{lin}$ and thus the conditions change as well. This implies that there is more flat outputs for
$\S_{aff}$ than for the associated $\S_{lin}$. Actually, the condition (FO1) applied to $\S_{lin}$ implies that
$(L_g\vp_0,L_g\vp_1)(x^*)\neq (0,0)$ {(thus obtaining the same necessary and sufficient conditions as those given in \cite{li2010flat} for two-input control-linear systems
 }), whereas (FO1) applied to $\S_{aff}$ still admits systems for which  $(L_g\vp_0,L_g\vp_1)(x^*)= (0,0)$ as the following example shows.

{\begin{example}  Consider the control-affine system:
\end{example}}
\vspace{-0.75cm}
 $$\begin{array}{lclcl}
  \dot z_0 = v_0 & & \dot z_1 &=& z_0 + z_2v_0\\
   & &\dot z_2 &=& z_3v_0\\
    & & &\vdots&\\
     & &\dot z_{k-1} &=& z_kv_0\\  
   & &\dot z_k & =& v_1
 \end{array}
$$
which is in the triangular form compatible with the chained form $TCh_1^k $. 
We claim that it is $x$-flat with $(\vp_0,\vp_1) =( z_1-z_0z_2, z_2)$ {as $x$-flat output around} $z^*=0$, although $(L_{g}\vp_0,L_{g}\vp_1)(0) = (0,0)$, for any vector field in $ \mathcal{G}$ such that $g(z^*)\not \in \mathcal{C}^{k-2}(z^*)$, provided that $v_0^* \neq 0 $ and $(1-z_3^*v_0^*)\neq 0$, the latter condition being always satisfied at $z^*=0$, but not in a neighborhood.

Indeed, we have 
$\dot \vp_0 = z_0 - z_0 z_3v_0$, $\dot \vp_1 = z_3v_0$ and it follows that $\dot \vp_0 = z_0(1-\dot \vp_1)$, from which we deduce $z_0 = \frac{\dot \vp_0}{1-\dot \vp_1}$, provided that $1-\dot \vp_1 = 1-z_3^*v_0^*\neq 0$. By differentiating that relation, we get $v_0 = \dot z_0 = \frac{d}{dt}(\frac{\dot \vp_0}{1-\dot \vp_1}) = \delta_0(\bar \vp_0^2,\bar \vp_1^2)$, where $\bar \vp_i^j = (\vp_i, \dot \vp_i, \cdots,\vp_i^{(j)} )$. From $\dot \vp_1 = z_3v_0$, we compute $z_3 = \frac{\dot \vp_1 }{v_0} = \gamma_3(\bar \vp_0^2,\bar \vp_1^2)$. Then, $\dot z_3$ gives $z_4 =  \gamma_4(\bar \vp_0^3,\bar \vp_1^3 )$ and so on. Finally we get $z_k =  \gamma_k(\bar \vp_0^{k-1},\bar \vp_1^{k-1} )$ and $v_1 = \delta_1(\bar \vp_0^k,\bar \vp_1^k )$. Thus $(\vp_0,\vp_1) =( z_1-z_0z_2, z_2)$ is indeed an $x$-flat output of the system around $z^*=0$ such that $z_3^*v_0^*\neq 1 $.

Let us now consider the chained form $Ch_1^k $ and take $g=g_0$. We compute $L_{g}\vp_0 = -z_0z_3v_0$, $L_{g}\vp_1 = z_3v_0 $ and we clearly have $(L_{g}\vp_0,L_{g}\vp_1)(0) = (0,0)$. Since the condition $(L_{g}\vp_0,L_{g}\vp_1)(z^*) \neq (0,0)$ is necessary for $(\vp_0,\vp_1)$ to be an $x$-flat output for  the chained form, see \cite{li2010flat},  we deduce that
$(\vp_0,\vp_1) =( z_1-z_0z_2, z_2)$ is not an $x$-flat {output} at $z^*=0$ for $Ch_1^k $. 
{$
   \hfill \qquad\square
$}

\bigskip

For control-linear systems $\S_{lin}$, the choice of a flat output is not unique (different choices are parameterized by an arbitrary function of three variables whose differentials annihilate $\mathcal{C}^{k-2}$, as assures Proposition \ref{chap 3.1 prop:output-eq:m=1} below) but all flat outputs exhibit the same singularity in control space {(see item $(F4)$ of Proposition \ref{chap 3.1 prop output m=1 S_lin})}, which is the control $u_{c}$, where $u_{c}\in U_{char}$ such that $u_{c,0}g_{0}+u_{c,1}g_{1} \in \mathcal{C}^{1}$ ( for any $x \in X$, it defines a one-dimensional linear subspace of $U=\mathbb{R}^{2}$). 
In the control-affine case, the nature of singularities changes substantially: each choice of a flat output creates its own singularities in the control space. More precisely, a flat output $(\varphi_{0},\varphi_{1})$ ceases to be a flat output for controls $u^{\ast}$ belonging to $U_{_{\mathcal{L}-sing}}$ which is the union of $\bigcup_{i=0}^{k-3}U_{sing}^{i}$ (universal for all choices of $(\varphi_{0},\varphi_{1})$ and consisting, for each fixed $x\in X$, of {the union of } $k-2$ one-dimensional affine subspaces of $U=\mathbb{R}^{2}$) and of $U_{_{\mathcal{L}-sing}}^{k-2}$, which is a one-dimensional affine subspace of $U=\mathbb{R}^{2}$ that depends on $(\varphi_{0},\varphi_{1})$ since $\mathcal{L} = (\mspan \{d\varphi_0,d\varphi_1\})^\perp$. 
All those $k-1$ affine subspaces are, in general, different although some of them may coincide and, indeed, in the control-linear case all of them coincide and reduce to the linear-space of $U=\mathbb{R}^{2}$ containing the characteristic controls $u_{c}$ that correspond {to }the characteristic distribution~$\mathcal{C}^{1}$, 
that is, the corresponding trajectories remain tangent to $\mathcal{C}^{1}$. Moreover, if we apply an invertible feedback $u=\b \t u$ (which always exists and can be explicitly calculated) such that $\mathcal{C}^{1} = \mspan\{\t g_1\}$ and $\mathcal{G}^{0} = \mspan\{\t g_0,\t g_1\}$, a control $\t u_c$ is characteristic, that is, singular for flatness of $\S_{lin}$, if and only if the feedback modified control is $\t u_c =  \b^{(-1)} u_c =(0,\t u_{c,1})^T$.

Now it is clear that {the control-affine system $\S_{aff}$} is flat if we avoid the universal singular set $\bigcup_{i=0}^{k-3}U_{sing}^{i}$ as well as  the set singular for all choices of flat outputs $(\varphi_{0},\varphi_{1})$, that is the set $\bigcap U_{_{\mathcal{L}-sing}}^{k-2}$ (the intersection taken over all $\mathcal{L}$), which explains different statements for a fixed choice of $(\varphi_{0},\varphi_{1})$ in item ($F2$)($ii$) and an arbitrary choice of $(\varphi_{0},\varphi_{1})$ in item ($F1$).

\medskip

Notice that Theorem \ref{chap 3.1 thm:output:m=1} is valid for any $k \geq 3$ (thus for a system defined on a manifold $X$ of dimension at least 4). In fact, in item {$(ii)$}, we use the characteristic distribution $\mathcal{C}^{k-2}$ of $\mathcal{G}^{k-2}$, but if dim $X$ = 3, i.e., $k=2$,  such a distribution does not exist and item {$(ii)$} does not apply to that case. Item {$(iii)$}, however, is well defined even for dim $X$ = 3 and remains equivalent to $(i)$ {.}
\medskip 

As an immediate corollary of Theorem \ref{chap 3.1 thm:output:m=1}, we obtain a system of first order PDE's, described by Proposition \ref{chap 3.1 prop:output-eq:m=1} below, whose solutions give all $x$-flat outputs.
Like for systems equivalent to the chained form (see \cite{li2010flat}), $x$-flat outputs for the systems feedback equivalent to the triangular form $TCh_1^k$ are far from being unique: since the distribution $\mathcal C^{k-2}$ is involutive and of corank three, there are as many functions $\varphi_0$ satisfying $L_c \varphi_0 =  0$, for any  $c \in \mathcal C^{k-2}$, as functions of three variables. 
{Indeed, according to the following proposition}, $\varphi_0$ can be chosen as any function of {the} three independent functions, whose differentials annihilate $\mathcal C^{k-2}$, and {if moreover, $<d\vp_0, \mathcal{G}^0>(x^*) \neq 0$},  then there exists a unique $\varphi_1$ (up to a diffeomorphism) completing it to an $x$-flat output.

\begin{proposition}\label{chap 3.1 prop:output-eq:m=1}
 Consider a two-input control-affine system $\Sigma_{aff}: \dot x = f(x) + u_0g_0(x) + u_1g_1(x)$, defined on a manifold $X$, of dimension $k+1\geq 4$,  that is locally, around $x^* \in X$, static feedback equivalent to $TCh_1^k$. Let $\mathcal{C}^{k-2} = \mspan\{c_1, \cdots, c_{k-2}\}$ be the characteristic distribution of $\mathcal{G}^{k-2}$ such that $c_{k-2}(x^*)\not\in \mathcal{C}^{k-3}(x^*)$ and $g$ be an arbitrary vector field in $\mathcal{G}$ such that $g(x^*)\not\in \mathcal{C}^{k-2}(x^*)$. Then
 \begin{enumerate}[(i)]
  \item For any smooth function $\vp_0$ such that
  $$
  \mbox{(Flat 1)} \quad L_{c_i}\vp_0=0, \, \dleq 1 i k-2, \mbox{ and } <d\vp_0, \mathcal{G}^{k-2}>(x^*) \neq 0,
  $$
  the distribution $\mathcal{L} =\mathcal{C}^{k-2} + \mspan\{v \}$ is involutive, where $v = (L_{g}\vp_0)[c_{k-2},g]-(L_{[c_{k-2},g]}\vp_0)g$.

  \item A pair $(\varphi_0, \varphi_1)$ of smooth functions defined on a neighborhood of $x^*$ is an $x$-flat output at $(x^*,u^*)$ with $u^* \not\in U_{\mathcal{L}-sing}(x^*)$, if and only if (after permuting $\varphi_0$ and $\varphi_1$, if necessary) $\varphi_0$ is any function satisfying
  (Flat 1) and $\varphi_1 $ satisfies
  $$\mbox{(Flat 2)} \quad \left\{
  \begin{array}{l@{}c@{}l}   
  (d\varphi_0 &\wedge& d\varphi_1 \wedge d\dot \varphi_0 \wedge d\dot \varphi_1)(x^*, u^*) \neq 0,\\
  L_{c_i}\vp_1&=&0, {\mbox{ for }} \dleq 1 i k-2, \\
  L_{v}\vp_1&= &0.
  \end{array}
  \right.
    $$
  \item If in   (Flat 1), we replace $<d\vp_0, \mathcal{G}^{k-2}>(x^*) \neq 0$ by $<d\vp_0, \mathcal{G}^{0}>(x^*) \neq 0$, then for any function $\varphi_0$ satisfying $L_c \varphi_0 =  0$, for any  $c \in \mathcal C^{k-2}$, and $<d\vp_0, \mathcal{G}^{0}>(x^*) \neq 0$, there always exists $\varphi_1$ such that the pair $(\varphi_0,\varphi_1)$ is an $x$-flat output of $\S_{aff}$; given any such $\varphi_0$, the choice of $\varphi_1$ is unique, up to a diffeomorphism, that is, if $(\varphi_0,\tilde \varphi_1)$ is another minimal $x$-flat output, then there exists a smooth map $h$, smoothly invertible with respect to the second argument, such that
$$
\tilde \varphi_1 =  h(\varphi_0,\varphi_1).
$$
 \end{enumerate}
\end{proposition}
\textit{Remark.} Notice that for a function $\vp_0$ satisfying  $<d\vp_0, \mathcal{G}^{k-2}>(x^*) \neq 0$ (and not the stronger condition $<d\vp_0, \mathcal{G}^{0}>(x^*) \neq 0$, or equivalently $L_{g}\vp_0(x^*) \neq 0 $, see Proposition \ref{chap 3.1 prop:output-eq:m=1}(iii)), it can be impossible to find, among all solutions of $L_{c_i}\vp_1=L_{v}\vp_1= 0$, $\dleq 1 i k-2$, a function $\vp_1$ satisfying $ (d\varphi_0 \wedge d\varphi_1 \wedge d\dot \varphi_0 \wedge d\dot \varphi_1)(x^*, u^*) \neq 0$ and therefore item $(iii)$ does not hold, in general, under the weaker condition $<d\vp_0, \mathcal{G}^{k-2}>(x^*) \neq 0$. This is, for example, the case of control-linear systems. 

As expected, the system of PDE's allowing us to compute all $x$-flat outputs of a system locally static feedback equivalent to $TCh_1^k$ does not depend on the drift $f$ and it is the same as that provided in \cite{li2010flat} for $x$-flat outputs in the case of control-linear $\Sigma_{lin}$ feedback equivalent to the chained form.  {For more details and the proof of Proposition~\ref{chap 3.1 prop:output-eq:m=1} in the case  $L_{g}\vp_0(x^*) \neq 0$, we refer the reader to \cite{li2010flat}.}

\bigskip

Finally, it turns out that almost all $x$-flat outputs are compatible with the triangular form $TCh_1^k$  (as are  $x$-flat outputs of the chained form). In fact, for any given flat output $(\vp_0,\vp_1)$ of a system $\Sigma_{aff}$ feedback equivalent to $TCh_1^k$,  {verifying $(L_{g}\vp_0,L_{g}\vp_1)(x^*) \neq (0,0)$}, we can bring $\Sigma_{aff}$ into $TCh_1^k$ for which $\vp_0$ and $\vp_1$ serve as the two top variables, as the following proposition assures.  The following result is technical and will be useful in our proofs, but it has its own interest.

\begin{proposition}\label{chap 3.1 prop-charac-FO}
Assume that $\Sigma_{aff}$ is locally, around $x^{\ast}$, static feedback equivalent to the triangular form $TCh_1^k$ and let $(\vp_0,\vp_1)$ be an $x$-flat output around $(x^{\ast}, u^{\ast})$, {such that $(L_{g}\vp_0,L_{g}\vp_1)(x^*) \neq (0,0)$, where $g$ is an arbitrary vector field in $\mathcal{G}$ such that $g(x^*)\not\in \mathcal{C}^{k-2}(x^*)$. Then we can bring $\Sigma_{aff}$ to $TCh_1^k$ around $z^{\ast}$ such that $z_{0}=\varphi_{0}$ and $z_{1}=\varphi_{1}$} (after permuting $\varphi_{0}$ and $\varphi_{1}$, if necessary).
\end{proposition}

{\textit{Remark.} The above proposition is valid around $z^{\ast}$ which is not necessary equal to~0. If we want to map $x^*$ into  $z^{\ast}=0$, then an affine transformation of flat outputs may be needed. More precisely,  we can bring $\Sigma_{aff}$ to $TCh_1^k$ around $z^{\ast}=0$ such that $z_{0}=\varphi_{0}$ and $z_{1}=\varphi_{1}+k_{0}\varphi_{0}$ (after permuting $\varphi_{0}$ and $\varphi_{1}$), where $k_{0}\in \mathbb{R}$.} 

{\subsection{Flatness of control  systems static feedback equivalent to $TCh_m^k$}}

We now turn to the case $m\geq 2$.
It is clear that $TCh_m^k$ is $x$-flat, with $\varphi = (z_0,z_1^1, \cdots,z_m^1 )$ being a flat output, at any point  $(z^*,v^*) \in \mathbb{R}^{km+1}\times \mathbb{R}^{m+1} $ satisfying
$$\label{chap 3.1 eq:regularity:m>=2}
\mrk F^l(z^*,v^*) = m, \mbox{ for } 1 \leq l \leq k-1,
$$
where $F^l=(F_{ij}^{l})$, for $ 1 \leq l \leq k-1,$ is the $m\times m$ matrix given by
$$
F^l_{ij} =\pfrac{(f^l_j + z^{l+1}_jv_{0})}{z^{l+1}_i} ,  \mbox{ for } 1 \leq i, j \leq m.
$$
Therefore, flat systems equivalent to $TCh_m^k$  exhibit singularities in the control space (depending on the state)
defined in an invariant way by 
$$
U_{m-sing}(x)=\bigcup_{i=0}^{k-2}U^i_{m-sing}(x),
$$
where
$$
U^i_{m-sing}(x)=\{u(x) \in\R^2:\mrk(\mathcal{G}^i+ [f+gu,\mathcal{L}^{i+1} ]) {(x)}< (i+2)m+1\},
$$
with $\mathcal{L}^{i+1} = \mathcal{C}^{i+1} $, 
for $\dleq 0 i k-3$, where $ \mathcal{C}^{i+1} $ is the characteristic distribution of $ \mathcal{G}^{i+1} $, and  $\mathcal{L}^{k-1} = \mathcal{L}$, the involutive {subdistribution} of  $ \mathcal{G}^{k-2} $
and $gu = \sum_{i=0}^m u_ig_i$. This singularity is excluded by 
item \textit{(m-F1)} of the next theorem describing all minimal $x$-flat outputs of control-affine systems feedback equivalent to the triangular form $TCh_m^k$.

\begin{theorem}
\label{chap 3.1 thm:flat_outputs:m>=2}
Consider a control-affine system $\Sigma_{aff}: \dot x = f(x) + \sum_{i=0}^m u_ig_i(x),$ with $m\geq~2$, defined on an open subset $X$ of $\mathbb{R}^{km+1}$, where $k\geq 2$,  that is locally, around $x^* \in~X$, static feedback equivalent to $TCh_m^k$ and its associated control-linear system $\S_{lin}: \dot x= \sum_{i=0}^m u_ig_i(x)$.
\begin{enumerate}[(m-F1)]
 \item $\Sigma_{aff}$ is $x$-flat, of differential weight $(k+1)(m+1)$, at any $(x^*,u^*)\in X\times \mathbb{R}^{m+1}$ such that $u^* \not\in  {U_{m-sing}}(x^*).$

 \item If $(\vp_0, \cdots, \vp_m)$ is a minimal $x$-flat output of $\Sigma_{aff}$ at $(x^*,u^*)$, where $u^*\not\in U_{m-sing}(x^*)$, then there exists an open neighborhood $\mathcal{X}^*$ of $x^*$ and coordinates
$(z_0,z_1^1,$ $ \cdots ,z_m^1, \cdots,$ $z_1^k, \cdots ,z_m^k)$ on  $\mathcal{X}^*$ in which $\S_{aff}$ is locally feedback equivalent to the triangular form  $TCh_m^k$, such that $\vp_0 = z_0$ and $\vp_i = z_i^1$, for $\dleq 1 i m$, after permuting the components $\vp_i$ of the flat output $\vp$, if necessary.

 \item Let $\vp_0$, $\vp_1, \cdots, \vp_m$ be $m+1$ smooth functions defined in a {neighborhood of $x^*$}. The following conditions are equivalent:
\begin{enumerate}[{(}i{)}]
\item The $(m+1)$-tuple $(\vp_0,\vp_1, \cdots, \vp_m)$ is a minimal $x$-flat output of $\Sigma_{aff}$ at $(x^*,u^*)$, where $ u^* \not\in U_{m-sing}(x^*)$;
\item The  $(m+1)$-tuple $(\vp_0,\vp_1, \cdots, \vp_m)$ is a minimal $x$-flat output of $\S_{lin}$ at $(x^*,\tilde u^*)$, where $\tilde u^*$ is such that $\sum_{i=0}^m \tilde u^*_{i}g_i(x^*)\not\in\mathcal{C}^{1}(x^*)$, where $\mathcal{C}^{1}$ is the characteristic distribution of $\mathcal{G}^1$;
\item The  $(m+1)$-tuple $(\vp_0,\vp_1, \cdots, \vp_m)$  satisfies the following conditions {in a neighborhood of $x^*$}:
    \begin{enumerate}[{(m-FO}1{)}]
		\item $d\varphi_0 \wedge d\varphi_1\wedge \cdots   \wedge d\varphi_m (x^*) \neq 0$;
		\item $\mathcal{L} = (\mspan \{d\varphi_0,d\varphi_1,\cdots, d\varphi_m \})^\perp$, where $\mathcal{L}$ denotes the involutive subdistribution of corank one in $\mathcal{G}^{k-2}$.
    \end{enumerate}
\end{enumerate}
Moreover, the  $(m+1)$-tuple $(\vp_0,\vp_1, \cdots, \vp_m)$ is unique, up to a diffeomorphism, i.e., if $(\tilde \varphi_0,\tilde \varphi_1, \cdots,\tilde \vp_m)$ is another minimal $x$-flat output, then there exist smooth maps $h_i$ such that $\t \vp_i = h_i(\vp_0,\vp_1, \cdots, \vp_m)$, $\dleq  0 i m$, and $h = (h_0,h_1, \cdots, h_m)$ is a {local} diffeomorphism.
\end{enumerate}
\end{theorem}

Theorem \ref{chap 3.1 thm:flat_outputs:m>=2} indicates how flatness of control-affine systems locally equivalent to $TCh_m^k$ reminds, but also how it differs from, that of control-linear systems locally equivalent to the $m$-chained form  {$Ch_m^k$}.

While Theorem \ref{chap 3.1 thm:output:m=1}, associated to the case $m=1$, allows us to compute all $x$-flat outputs of $TCh_1^k$, Theorem \ref{chap 3.1 thm:flat_outputs:m>=2} describes all minimal $x$-flat outputs of $TCh_m^k$. Functions whose differentials annihilate~$\mathcal{L}$  are clearly not the only $x$-flat outputs of $TCh_m^k$. They are, however, the only that possess the minimality property, i.e., when determining, with their help, all state and control variables, we use the minimal possible number of derivatives, which is $(k+1)(m+1)$, {see the proof of Theorem \ref{chap 3.1 thm:flat_outputs:m>=2}}. 
According to item $(ii)$, their description coincides with that of minimal $x$-flat outputs of $\S_{lin}$. Indeed, conditions \textit{(m-FO1)-(m-FO2)} are the same as those given in \cite{respondek2003symmetries} for control-linear systems feedback equivalent to the $m$-chained form. The presence of the drift has no influence on characterizing minimal $x$-flat outputs, but, analogously to the case $m=1$, it plays a role in describing singularities in the control space.

As for the characterization of the $m$-chained form and, consequently, of control-affine systems static feedback equivalent to $TCh_m^k$, the involutive subdistribution $\mathcal{L}$ of corank one in $\mathcal{G}^{k-2}$ is crucial for minimal $x$-flat outputs computation. Indeed, all minimal $x$-flat outputs are determined by~$\mathcal{L}$.
In contrast with the case $m=1$, where the choice of $x$-flat outputs is parameterized by a function of three well chosen variables, minimal $x$-flat outputs of $TCh_m^k$ are unique (as they are for the $m$-chained form). This is a consequence of the uniqueness of the involutive subdistribution~$\mathcal{L}$ of corank one in $\mathcal{G}^{k-2}$, in the case $m\geq 2$, and multiple noncanonical choices of $\mathcal{L}$, if $m=1$.

For control-affine systems, it is  the drift $f$, the characteristic distributions $\mathcal{C}^i$, for $\dleq 1 i k-2$, and  the involutive subdistribution  $\mathcal{L}$ of corank one in $\mathcal{G}^{k-2}$, that describe singularities in the control space. Although $\mathcal{L}$ is not involved in the compatibility conditions (see item \textit{(m-Comp)} of Theorem~\ref{chap 3.1 thm:m>=2}), it plays an important role in determining the singular controls at which the system ceases to be flat.

The description of the set of singular controls $U_{m-sing}$ is also valid  for driftless systems, i.e., for $f=0$, but it is  redundant. In fact, the set of singular controls $u_{c}$ for control-linear systems can be described using the first characteristic distribution $\mathcal{C}^1$ only: the singular controls $u_{c}$ are such that the corresponding trajectories are tangent to the characteristic distribution $\mathcal{C}^{1}$, that is, $u_{c}$ verifying $\sum_{i=0}^m u_{c,i}(x)g_i(x) \in\mathcal{C}^{1}(x)$. Clearly, they form, for any $x\in X$, an $m$-dimensional linear subspace of $U=\mathbb{R}^{m+1}$. If we apply an invertible feedback $u=\beta\tilde{u}$ such that $\mathcal{C}^{1}=\mspan\{\t g_1, \cdots, \t g_m\}$ and $\mathcal{G}^0 = \mspan\{\t g_0\} + \mathcal{C}^{1}$, then the singular controls $\tilde{u}_{c}$ are of the form $\tilde{u}_{c} = (0,\tilde{u}_{c,1},\cdots, \tilde{u}_{c,m})$.
\medskip

Finally, it turns out that minimal $x$-flat outputs and the triangular form $TCh_m^k$ are compatible: in fact, for any $m+1$ smooth functions $\vp_0,\vp_1, \cdots, \vp_m$ that form a minimal $x$-flat output of a system $\S_{aff}$ feedback equivalent to $TCh_m^k$, we can bring $\S_{aff}$ into the form $TCh_m^k$ for which $\vp_0,\vp_1, \cdots, \vp_m$ play the role of the top variables, {as item \textit{(m-F2)} assures. An analogous result}  is also valid for minimal $x$-flat outputs and the $m$-chained form, see \cite{li2011geometry}.
\medskip

As an immediate corollary of Theorem \ref{chap 3.1 thm:flat_outputs:m>=2}, we get the following system of PDE's whose solutions give all minimal $x$-flat outputs for control-affine systems static feedback equivalent to $TCh_m^k$.
Denote by $v_j$, for $\dleq 1 j (k-1)m$, the vector fields spanning the distribution $\mathcal{L}$ (for their computation see Appendix A). 

\begin{proposition}
\label{chap 3.1 prop:flat_outputs_eq:m>=2}
Consider a control-affine system $\Sigma_{aff}: \dot x = f(x) + \sum_{i=0}^m u_ig_i(x),$ with $m\geq 2$, defined on an open subset $X$ of $\mathbb{R}^{km+1}$, where $k\geq 2$, that is locally, around $x^* \in~X$, static feedback equivalent to $TCh_m^k$. Let $\mathcal{L} = \mspan \{v_j$, \, $\dleq 1 j (k-1)m\}$ be the involutive subdistribution of corank one in $\mathcal{G}^{k-2}$. Then smooth functions $\vp_0$, $\vp_1, \cdots, \vp_m$, defined in a neighborhood of $x^*$, form a minimal $x$-flat output at $(x^*,u^{\ast})$, $u^{\ast} \not\in U_{m-sing}(x^{\ast})$ if and only if
$$L_{v_j}\vp_i=  0, \; \dleq 1 j (k-1)m, \; \dleq 0 i m, $$
and $d\varphi_0 \wedge d\varphi_1\wedge \cdots   \wedge d\varphi_m (x^*) \neq 0$.
\end{proposition}

\section{Examples and applications}\label{chap 3.1 sec:application}
\subsection{Example: $TCh_{1}^k$ is not necessary for flatness}\label{chap 3.1 ssec_academical example}

In the previous section we have seen that systems locally static feedback equivalent to the triangular form {to $TCh_m^k$, $m=1$ or $m\geq 2$}, are $x$-flat and we have described all $x$-flat outputs.
Therefore being static feedback equivalent {to $TCh_m^k$, $m=1$ or $m\geq 2$} is sufficient for $x$-flatness. A natural question arises: is static feedback equivalence {to $TCh_m^k$} necessary for flatness, provided that the control-linear subsystem is static feedback equivalent to the chained form? The next example gives a negative answer to this question.  Consider the following control-affine system whose control-linear part is already in the chained form {$Ch_1^4$, but whose drift $f$ does not satisfy the compatibility condition $(Comp)$ and thus the system cannot be transformed into $TCh_{1}^4$}:
$$
\left\{
\begin{array}{l lcc c l}
 \dot z_0 = v_0 & \dot z_1 &=& z_3 &+& z_2v_0\\
                & \dot z_2 &=& -z_4 &+& z_3v_0\\
                & \dot z_3 &=& a(\bar z_{3}) &+& z_4v_0\\
                & \dot z_4 &=& v_1
\end{array}
\right.
$$
where $a$ is a smooth function depending on $z_0, z_1 , z_2, z_{3}$. The pair $(\varphi_0, \varphi_1) = (z_0, z_1)$ is an $x$-flat output. Indeed, we have $\varphi_0 = z_0$ implying $\dot{\varphi}_0 = v_0$ and $\varphi_1 = z_1$ implying
$$
\begin{array}{l}
 \dot \varphi_1 = z_3 + z_2v_0  = z_3 + z_2 \dot \varphi_0\\
 \ddot \varphi_1 =a(\varphi_0,\varphi_1, z_2,z_{3})+ z_3 \dot \varphi_0^2+ z_2\ddot \varphi_0.
 \end{array}
$$
These expressions allow us to calculate $z_2$ and $z_3$ via the implicit function theorem as
$$
\begin{array}{l}
z_2 = \gamma_2(\bar \varphi_0^2, \bar \varphi_1^2)\\[2mm]
z_3 = \gamma_3(\bar \varphi_0^2, \bar \varphi_1^2),
 \end{array}
$$
for some functions $\gamma_2$, $\gamma_3$, where $\bar \varphi^l$ denotes $(\varphi, \dot \varphi, \cdots, \varphi^{(l)})$. By differentiating $z_3$, we deduce $z_4 = \gamma_4(\bar \varphi_0^3, \bar \varphi_1^3)$ which yields $ v_1 = \delta_1 (\bar \varphi_0^4, \bar \varphi_1^4)$. So we have expressed all state and control variables as functions of $\varphi_0$ and $\varphi_1$ and their derivatives proving that $(\varphi_0, \varphi_1) = (z_0, z_1)$ is, indeed, an $x$-flat output.

\subsection{Application to mechanical systems: coin rolling without slipping on a moving table}

\begin{figure}[htbp]
\centering
\includegraphics[width= 6.25cm]{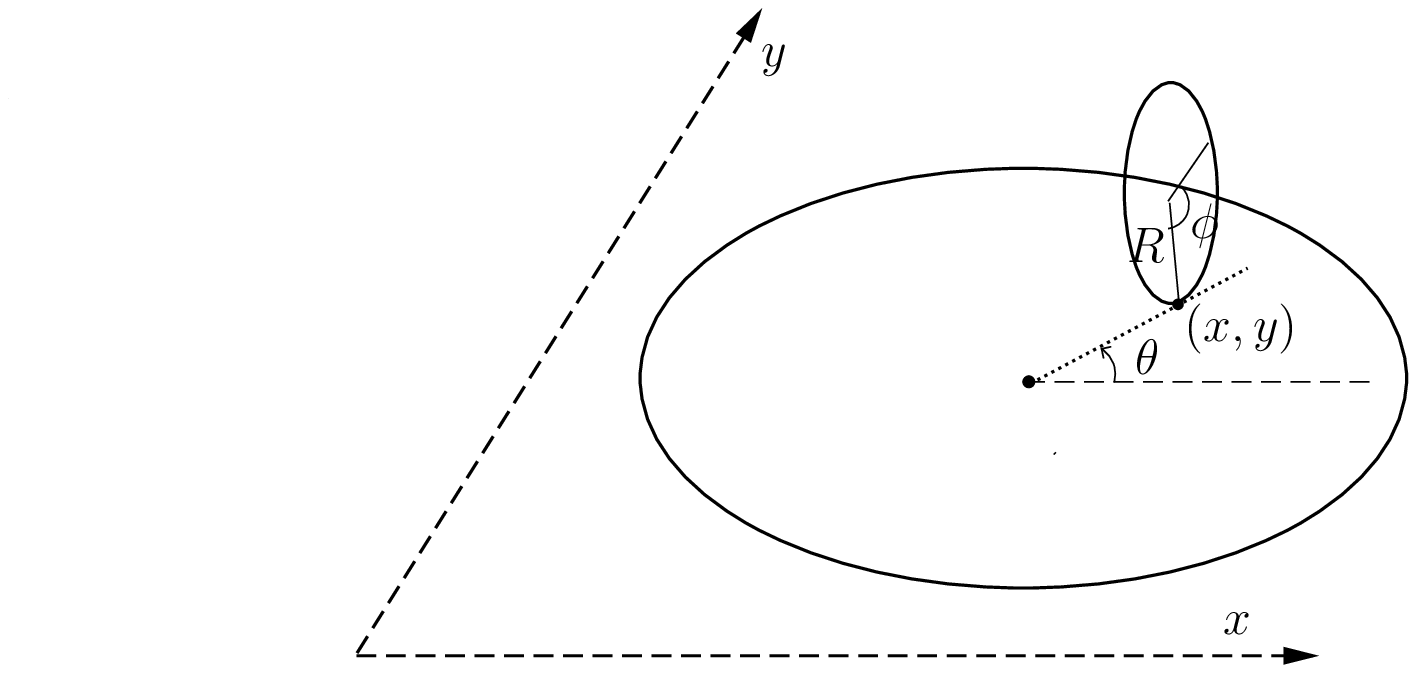}
\caption{The coin on a moving table}\label{chap 3.1 Fig-1}
\end{figure}

Consider a vertical coin of radius $R$ rolling without slipping on a moving table, see Figure~\ref{chap 3.1 Fig-1}. Assume that the surface of the table is on the $xy$-plane and denote by $(x,y)$ the position of the contact point of the coin with the table, and by $\theta$ and $\phi$, respectively, the orientation of the {vertical plane containing the} coin and the rotation angle of the coin. Then the configuration space for the system is $Q=SE(2) \times S^{1}$  and is parameterized by the generalized coordinates $q=((x,y,\theta),\phi)$.

Assume that the table {moves with respect to the inertial frame} obeying the differential equations
\begin{equation}\label{chap 3.1 equ-table}
\left.
     \begin{array}{lll}
       \dot{x}{_t} &=&  \alpha(x{_t},y{_t})  \\
       \dot{y}{_t} &=&  \beta(x{_t},y{_t}).
       \end{array}
  \right.
\end{equation}
for a smooth vector field $(\alpha,\beta)^{\top}$ on $\mathbb{R}^{2}$.

Therefore the nonholonomic constraints of rolling without slipping can be represented by
\begin{equation}\label{chap 3.1 constraint}
     \begin{array}{cll}
       \dot{x}\sin \theta  - \dot{y} \cos\theta &=&  0  \\
      (\dot{x}-\alpha)\cos \theta  + (\dot{y}-\beta) \sin\theta  &=& R \dot{\phi},
     \end{array}
\end{equation}
which leads to the kinematic model of the coin on a moving table as
\begin{equation} \label{chap 3.1 g-coin}
\Sigma_{coin}:\ \
 \left(
      \begin{array}{c}
        \dot{x} \\
        \dot{y} \\
        \dot{\theta} \\
        \dot{\phi} \\
        \end{array}
     \right) = \left(
                      \begin{array}{c}
                      \cos \theta (\alpha \cos \theta + \beta\sin\theta)\\
                      \sin \theta (\alpha \cos \theta + \beta\sin\theta)\\
                      0\\
                      0\\
                      \end{array}
                \right)
                +\left(
                      \begin{array}{c}
                      0\\
                      0\\
                       1\\
                      0\\
                      \end{array}
                \right)
                u_1 + \left(
                        \begin{array}{c}
                          R \cos \theta\\
                          R \sin \theta \\
                          0 \\
                          1 \\
                          \end{array}
                          \right)u_2.\\
                      \end{equation}
The system is control-affine because the nonholonomic constraints are affine (and not linear) as a result of the motion of the table with respect to the inertial frame.
\begin{remark}
Assume that $\alpha=-\omega y{_t}$, $\beta = \omega x{_t}$, that is, the motion equation of the table is
\begin{equation*}
\left.
     \begin{array}{llc@{}l@{}}
       \dot{x}{_t} &=&  -&\omega y{_t} \\
       \dot{y}{_t} &=&  &\omega x{_t},
       \end{array}
  \right.
\end{equation*}
meaning that the table rotates around its center point with the angular velocity $\omega$. Substituting $\alpha=-\omega y$, $\beta = \omega x$ into (\ref{chap 3.1 g-coin}), we obtain the model of the coin on a rotating table as
\begin{equation} \label{chap 3.1 coin-rotating-table}
 \left(
      \begin{array}{c}
        \dot{x} \\
        \dot{y} \\
        \dot{\theta} \\
        \dot{\phi} \\
        \end{array}
     \right) = \left(
                      \begin{array}{c}
                      \omega \cos \theta (x\sin\theta - y \cos \theta)\\
                      \omega \sin \theta (x\sin\theta - y \cos \theta)\\
                      0\\
                      0\\
                      \end{array}
                \right)
                +\left(
                      \begin{array}{c}
                      0\\
                      0\\
                       1\\
                      0\\
                      \end{array}
                \right)
                u_1 + \left(
                        \begin{array}{c}
                          R \cos \theta\\
                          R \sin \theta \\
                          0 \\
                          1 \\
                          \end{array}
                          \right)u_2,
                      \end{equation}
which coincides with the model given by \cite{kai2006extended}.
\end{remark}

\begin{proposition}
The coin on a moving table $\Sigma_{coin}$, given by (\ref{chap 3.1 g-coin}), is feedback equivalent to the triangular form $TCh_1^3$ if and only if the motion of the table is described by
$$
\left\{
     \begin{array}{llc}
       \dot{x}{_t} &=&  cy{_t}+ d  \\
       \dot{y}{_t} &=&  -cx{_t} + e
       \end{array}
  \right.
$$
where $c,$ $d,$ $e\in \mathbb{R}$ are constant.
\end{proposition}
\begin{remark}
Notice that introducing $\tilde{x}{_t}=x{_t}- e/c$ and $\tilde{y}{_t}=y{_t}+ d/c$, we obtain: 
\begin{equation*}
\left.
     \begin{array}{lll}
       \dot{\tilde{x}}{_t} &=&  c\tilde{y}{_t}\\
       \dot{\tilde{y}}{_t} &=&  -c\tilde{x}{_t}.
       \end{array}
  \right.
\end{equation*}
 The only motions of table that lead to the triangular form $TCh_1^3$ are thus constant speed rotations around a fixed point $( e/c,- d/c)$.
\end{remark}

\begin{proof}
{The system} $\Sigma_{coin}$ is feedback equivalent to the triangular form $TCh_1^k$ {if and only if} it satisfies the conditions \textit{(Ch1)-(Ch3)} and \textit{(Comp)} of Theorem~\ref{chap 3.1 thm:m>=2} or, equivalently, conditions \textit{(Ch1)'-(Ch2)'} and \textit{(Comp)}.
Consider the associated distribution {$\mathcal{G}$ and the drift $f$ given by:}
$$\mathcal{G} = \textrm{span} \,\{g_{1},g_{2}\}  = \textrm{span} \, \left \{
\left ( \begin{array}{c}
0\\
0\\
1 \\
0\\
\end{array} \right ),
\left ( \begin{array}{c}
 R \cos \theta\\
 R \sin \theta \\
 0 \\
 1 \\
\end{array} \right )
\right \}
{\mbox{ and } f =  \left ( \begin{array}{c}
 \cos \theta (\alpha \cos \theta + \beta\sin\theta)\\
\sin \theta (\alpha \cos \theta + \beta\sin\theta)\\
0\\
0\\
\end{array} \right ).}
$$
A straightforward calculation shows that
$$
g_{3}= [g_{1},g_{2}] = \left ( \begin{array}{c}
 -R \sin \theta\\
 R \cos \theta \\
 0 \\
 0 \\
\end{array} \right ),
 \ \
g_{4}= [g_{1},g_{3}] = \left ( \begin{array}{c}
 -R \cos \theta\\
 -R \sin \theta \\
  0 \\
 0 \\
\end{array} \right ).
$$
Therefore $\mathcal{G}^{1}=\mathcal{G}_{1}=\textrm{span} \,\{g_{1},g_{2},g_{3}\}$ and $\mathcal{G}^{2}=\mathcal{G}_{2}=\textrm{span} \,\{g_{1},g_{2},g_{3},g_{4}\}$ which gives that $\mrk \mathcal{G}^{1}=\mrk \mathcal{G}_{1}=3$ and $ \mrk\mathcal{G}^{2} =\mrk\mathcal{G}_{2}=4$ and thus conditions \textit{(Ch1)'-(Ch2)'} hold. Moreover, it is easy to see that $\mathcal{C}{^{1}}=\textrm{span}\, \{c\}$ where $c=g_{2}$ and a direct computation gives
$$
[f,c]=[f,g_{2}] = -\left ( \begin{array}{c}
 \gamma R \cos \theta \\
 \gamma R \sin \theta \\
 0 \\
 0 \\
\end{array} \right ),
$$
where
$$ \gamma =  \cos \theta \left ( \frac{\partial \alpha}{\partial x} \cos \theta + \frac{\partial \beta}{\partial x} \sin \theta \right) +  \sin \theta \left ( \frac{\partial \alpha}{\partial y} \cos \theta + \frac{\partial \beta}{\partial y} \sin \theta \right).
$$
The condition $(Comp)$ of Theorem~\ref{chap 3.1 thm:m>=2} requires that $[f,c] \subset \mathcal{G}^{1}$ implying that the vector fields $[f,c]$ and $g_{3}$ are colinear and this is the case if and only if $\gamma \equiv 0$. We thus have to solve
\begin{equation*}
 \cos \theta \left ( \frac{\partial \alpha}{\partial x} \cos \theta + \frac{\partial \beta}{\partial x} \sin \theta \right) + \sin \theta \left ( \frac{\partial \alpha}{\partial y} \cos \theta + \frac{\partial \beta}{\partial y} \sin \theta \right) =0.
\end{equation*}
Dividing the above equation by $\cos^{2}\theta$ and denoting $w=\tan \theta$, we get
$$
\frac{\partial \alpha}{\partial x} + \left(\frac{\partial \alpha}{\partial y}+\frac{\partial \beta}{\partial x}\right) w + \frac{\partial \beta}{\partial y} w^{2} = 0,
$$
which implies that
$$
\frac{\partial \alpha}{\partial x}=0, \quad \frac{\partial \beta}{\partial y}=0, \quad \frac{\partial \alpha}{\partial y}=-\frac{\partial \beta}{\partial x}.
$$
We get
$\alpha = \alpha(y),$ $ \beta = \beta (x) $ and then by the equality $\frac{\partial \alpha}{\partial y}=-\frac{\partial \beta}{\partial x}$, we have
$$
\alpha^{\prime}(y)=-\beta^{\prime}(x)=c,
$$
where $c\in\mathbb{R}$ is a constant. This gives
\begin{equation*}
\left.
     \begin{array}{lll}
       \alpha &=&  cy+ d   \\
       \beta &=&  -cx + e
       \end{array}
  \right.
\end{equation*}
where $c,e,f \in\mathbb{R}$ are constants and the motion of the table is described by
\begin{equation}\label{chap 3.1 equ-moving}
\left.
     \begin{array}{lll}
       \dot{x}{_t} &=&  cy{_t}+ d  \\
       \dot{y}{_t} &=&  -cx{_t} + e,
       \end{array}
  \right.
\end{equation}
which proves the proposition.
\end{proof}

\section{Proofs}\label{chap 3.1 sec:proofs}

\bigskip
%
%

\subsection{Proof of Theorem \ref{chap 3.1 thm:m=1}}

\begin{proof}\textit{Necessity.} Consider a two-input control-affine system $\Sigma_{aff}: \dot x = f(x) +u_0 g_0(x) +u_1 g_1(x) $ locally, around $x^*$, static feedback equivalent to $TCh_1^k$ and bring it into the form $TCh_1^k$, around~$z^*$. By abuse of notation, we continue to denote by $f$, $g_0$ and $g_1$, the drift and the controlled vector fields of $TCh_1^k$. The distribution $\mathcal{G} = \mspan \{ g_0,g_1 \}$, associated to $TCh_1^k$, is given by
$$
\mathcal{G} = \mspan \{\pfrac{}{z_k}, \pfrac{}{z_0}+ z_2\pfrac{}{z_1}+ \cdots + z_k\pfrac{}{z_{k-1}} \}.
$$
By an induction argument, it is immediate to show that
$$
\mathcal{G}^i = \mathcal{G}_i = \mspan \{\pfrac{}{z_{k-i}}, \cdots, \pfrac{}{z_k}, \pfrac{}{z_0}+ z_2\pfrac{}{z_1}+ \cdots + z_{k-i}\pfrac{}{z_{k-i-1}} \}.
$$
Thus $\mathcal{G}^{k-1} = TX$ and the distribution $\mathcal{G}^{k-3} $ is of constant rank $k-1$.
The characteristic distribution~$\mathcal{C}^i$ of $\mathcal{G}^i$ is given by
$$
\mathcal{C}^i =  \mspan \{\pfrac{}{z_{k-i+1}}, \cdots, \pfrac{}{z_k}\}, \; \dleq 1 i k-2.
$$
So it is immediate to see that $\mathcal{C}^{k-2}$ is contained in $\mathcal{G}^{k-3} $, this inclusion is of corank one and $\mathcal{G}^0(z^*) \not\subset \mathcal{C}^{k-2}(z^*)$.  This shows \textit{(Ch1)-(Ch3)}.

Moreover, we have
$$
[\pfrac{}{z_{k}}, f] = \pfrac{f_{k-1}}{z_{k}}\pfrac{}{z_{k-1}} \in \mathcal{G}^1
$$
and
$$
[\pfrac{}{z_{k-i+1}}, f] = \pfrac{f_{k-i}}{z_{k-i+1}}\pfrac{}{z_{k-i}} \mmod \mspan \{\pfrac{}{z_{k-i+1}}, \cdots, \pfrac{}{z_{k}}\},
$$
which is clearly in $\mathcal{G}^i $, for any $\dleq 2 i k-2$. It follows that $[f, \mathcal{C}^i] \subset  \mathcal{G}^i$, for $\dleq 1 i k-2$, which shows \textit{(Comp)}.
The conditions $(Ch1)-(Ch3)$ involve the distribution $\mathcal{G}$ only, so they are invariant under feedback of the form {$g\rightarrow g\b$.}
Obviously, $[g_{j},\mathcal{C}^{i}]\in \mathcal{G}^{i}$ (since $\mathcal{C}^{i}$ is characteristic for $\mathcal{G}^{i}$), for $\dleq 0 j 1$, $1\leq i \leq k-2$, and thus $(Comp)$ is invariant under feedback of the form $f\mapsto f+\alpha_{0}g_{0}+\alpha_{1}g_{1}$.

\bigskip

\textit{Sufficiency. } Consider a two-input control-affine system $\Sigma_{aff}: \dot x = f(x) +u_0 g_0(x) +u_1 g_1(x) $ satisfying the conditions \textit{(Ch1)-(Ch3)} and \textit{(Comp)}. As proved in \cite{pasillas2001contact}, the items \textit{(Ch1)-(Ch3)} assure the  existence of an invertible static feedback transformation $u= \b \t u$ and a change of coordinates $z = \phi(x)$
bringing the distribution $\mathcal{G}^{0}$ into the chained form, which transform the system $\Sigma_{aff}$ into
$$
\left\{
\begin{array}{l lcl c l}
 \dot z_0 =a_0(z)+ \t u_0 & \dot z_1 &=& a_1(z) &+& z_2\t u_0\\
                          &    & \vdots &  & & \\
                          & \dot z_{k-1} &=& a_{k-1} (z) &+& z_k\t u_0\\
                          & \dot z_k &=& a_{k} (z) &+& \t u_1
\end{array}
\right.
$$
with $a_i$ smooth functions.
Applying the invertible static feedback $ v_0 = a_0(z)+ \t u_0$ and $v_1 = a_k(z)+ \t u_1$, we obtain
$$
\left\{
\begin{array}{l lcl c l}
 \dot z_0 =v_0 & \dot z_1 &=& f_1(z) &+& z_2v_0\\
               &   & \vdots & & &\\
               &  \dot z_{k-1} &=& f_{k-1} (z) &+& z_kv_0\\
               &  \dot z_k &=& v_1 & &
\end{array}
\right.
$$
where $f_i = a_i - z_{i+1}a_0$.
In these coordinates, we have
$$
\mathcal{G}^i = \mathcal{G}_i = \mspan \{\pfrac{}{z_{k-i}}, \cdots, \pfrac{}{z_k}, \pfrac{}{z_0}+ z_2\pfrac{}{z_1}+ \cdots + z_{k-i}\pfrac{}{z_{k-i-1}} \}, \; \dleq 0 i k-1,
$$
and
$$
\mathcal{C}^i =  \mspan \{\pfrac{}{z_{k-i+1}}, \cdots, \pfrac{}{z_k}\}, \; \dleq 1 i k-2.
$$
From  $[f, \mathcal{C}^i] \subset  \mathcal{G}^i$, for any $\dleq 1 i k-2$, it follows immediately that
$$
\pfrac{f_i}{z_j} = 0, \mbox{ for }\dleq{i+2}j k\mbox{  and } \dleq 1 i k-2,
$$
which gives the triangular normal form $TCh_1^k$.
\end{proof}

\subsection{Proof of Theorem \ref{chap 3.1 thm:m>=2}}
\begin{proof}\textit{ Necessity.}
Consider a control-affine system $\Sigma: \dot x = f(x) +\sum_{i=0}^m u_i g_i(x)$ locally, around $x^*$, static feedback equivalent to $TCh_m^k$ and bring it into the form $TCh_m^k$, around~$z^*$. To simplify the notation, we continue to write $f$ and $g_i$, $\dleq 0 i m$, for the drift and the controlled vector fields of $TCh_m^k$ and we denote
$$
\mspan \{\pfrac{}{z^i}\} = \mspan \{\pfrac{}{z^i_1}, \cdots, \pfrac{}{z^i_m}\}.
$$
The distribution $\mathcal{G}^{0} = \mspan \{ g_i, \; \dleq 0 i m \}$, associated to $TCh_m^k$, is given by
$$
\mathcal{G}^{0} = \mspan \{g_0, \pfrac{}{z^k} \}.
$$
By an induction argument, it is immediate that
$$
\mathcal{G}^{i} = \mspan \{\pfrac{}{z^{k-i}}, \cdots, \pfrac{}{z^k}, g_0 \}, \; \dleq 0 i k-1.
$$
It follows that  $\mathcal{G}^{k -1}= TX$, the distribution $\mathcal{G}^{k-2}$ has constant rank $(k-1)m+1$ and contains an involutive subdistribution of constant corank one given by
$$
\mathcal{L}=  \mspan \{\pfrac{}{z^{2}}, \cdots, \pfrac{}{z^k}\},
$$
and $\mathcal{G}^0(z^*)$ is not contained in $\mathcal{L}(z^*)$. This shows \textit{(m-Ch1)-(m-Ch3)}.
The characteristic distribution of $\mathcal{G}^i$ is given by
$$
\mathcal{C}^i =  \mspan \{\pfrac{}{z^{k-i+1}}, \cdots, \pfrac{}{z^k}\}, \; \dleq 1 i k-2,
$$
and we have, for any $k-i+1\leq l \leq k$ and $1\leq j\leq m$,
$$
[\pfrac{}{z^{l}_j}, f] = \pfrac{f^{l-1}_1}{z^{l}_j}\pfrac{}{z^{l-1}_1}+ \cdots +\pfrac{f^{l-1}_m}{z^{l}_j}\pfrac{}{z^{l-1}_m} \ \textrm{mod} \ \mathcal{C}^{i}
$$
which is clearly in $\mathcal{G}^{i}$. Thus
$[f, \mathcal{C}^i] \subset  \mathcal{G}^i$, for $1\leq i\leq k-2$, which proves item $(m-Comp)$.
\bigskip

\textit{Sufficiency.} Consider the control-affine system $\Sigma_{aff}: \dot x = f(x) +\sum_{i=0}^m u_i g_i(x)$ satisfying the conditions \textit{(m-Ch1)-(m-Ch3)} and \textit{(m-Comp)}.
According to Theorem 5.6 in \cite{pasillas2001contact}, the items \textit{(m-Ch1)-(m-Ch3)} assure the  existence of an invertible static feedback transformation $u= \b \t u$ and of a change of coordinates $z = \phi(x)$ (see Appendix B
where we explain how to construct the diffeomorphism $\phi$ and the feedback transformation) bringing the distribution~$\mathcal{G}^{0}$ into the $m$-chained form and thus the system $\Sigma_{aff}$ into
$$
\left\{
\begin{array}{l lcl c lcl}
 \dot z_0 =a_0(z)+\t u_0 & \dot z_1^1&=& a_1^1(z) + z_1^2\t u_0& \cdots & \dot z_m^1 &=& a_m^1( z) + z_m^2\t u_0\vspace{0.2cm}\\
                 & \dot z_1^2 &=& a_1^2( z) + z_1^3\t u_0&   & \dot z_m^2 &=& a_m^2( z) + z_m^3\t u_0\\
                 &            &\vdots&    & & & \vdots& \\
                 & \dot z_1^{k-1} &=& a_1^{k-1}( z) + z_1^k\t u_0& \cdots  & \dot z_m^{k-1} &=& a_m^{k-1}( z) + z_m^k\t u_0\vspace{0.2cm}\\
                 & \dot z_1^{k} &=& a_1^{k}( z) + \t u_1& \cdots  & \dot z_m^{k} &=& a_m^k( z) +\t u_m
\end{array}
\right.
$$
with $a^i_j$ smooth functions.
Applying the invertible static feedback $ v_0 = a_0(z)+ \t u_0$ and $ v_i = a_i^k(z)+ \t u_i$, for $\dleq 1 i m$, we get
$$
\left\{
\begin{array}{l lcl c lcl}
 \dot z_0 =v_0 & \dot z_1^1&=& f_1^1(z) + z_1^2v_0& \cdots & \dot z_m^1 &=& f_m^1(z) + z_m^2v_0\vspace{0.2cm}\\
                 & \dot z_1^2 &=& f_1^2(z) + z_1^3v_0&   & \dot z_m^2 &=& f_m^2(z) + z_m^3v_0\\
                 &            &\vdots&    & & & \vdots& \\
                 & \dot z_1^{k-1} &=& f_1^{k-1}(z) + z_1^kv_0&  \cdots & \dot z_m^{k-1} &=& f_m^{k-1}(z) + z_m^kv_0\vspace{0.2cm}\\
                 & \dot z_1^{k} &=& v_1& \cdots  & \dot z_m^{k} &=& v_m
\end{array}
\right.
$$
with $ f^i_j = a^i_j - z_j^{i+1}a_0$..
In the $z$-coordinates, we have
$$
\mathcal{G}^i = \mspan \{\pfrac{}{z^{k-i}}, \cdots, \pfrac{}{z^k}, g_0 \}, \; \dleq 0 i k-1.
$$
The characteristic distribution of $\mathcal{G}^i$ is given by
$$
\mathcal{C}^i =  \mspan \{\pfrac{}{z^{k-i+1}}, \cdots, \pfrac{}{z^k}\}, \; \dleq 1 i k-2,
$$
and the corank one involutive subdistribution of $\mathcal{G}^{k-2}$ by
$$
\mathcal{L}=  \mspan \{\pfrac{}{z^{2}}, \cdots, \pfrac{}{z^k}\}.
$$
We have, for $\dleq 1 i k-2$,
$$
[\pfrac{}{z^{k-i+1}_j}, f] =\sum_{l=1}^{m}\sum_{s=1}^{k-i-1} \pfrac{f^{s}_l}{z^{k-i+1}_j}\pfrac{}{z^{s}_l}  \mmod \mspan \{\pfrac{}{z^{k-i}}, \cdots, \pfrac{}{z^{k-1}}\}
$$
and since $[\pfrac{}{z^{k-i+1}_j}, f] \in \mathcal{G}^i$, for any $\dleq 1 j m$, we obtain
$$
\frac{f^{s}_l}{z^{k-i+1}_j}=0,\ \textrm{for any} \ \dleq 1 {j,l} m, \ \dleq 1 s k-i-1.
$$
It follows that $f$ exhibits the desired trangular form $TCh_m^k$.
\end{proof}

\subsection{Proof of Theorem \ref{chap 3.1 thm:output:m=1}}
\begin{proof}
[Proof of (F1)] Consider the two-input control-affine system $\Sigma: \dot x = f(x) +u_0 g_0(x) +u_1 g_1(x) $ locally, around $x^*$, feedback equivalent to $TCh_1^k$ and bring it into the form $TCh_1^k$, around $z^*$. To simplify notation, we continue to denote by $f$, respectively by $g_0$ and $g_1$, the drift, respectively the controlled vector fields of $TCh_1^k$.

It is clear that $TCh_1^k$ is $x$-flat, with $\varphi = (z_0,z_1)$ being a flat output, at any point  $(z^*,v^*)$ satisfying
$$
\pfrac{f_i}{z_{i+1}}(z^*)+v^*_{0}\neq 0,  \mbox{ for } 1 \leq i \leq k-1, 
$$
where $v^* = (v^*_0,v^*_1)$.
Recall that, in coordinates $z$, we have
$$
\mathcal{G}^i = \mspan \{\pfrac{}{z_{k-i}}, \cdots, \pfrac{}{z_k}, \pfrac{}{z_0}+ z_2\pfrac{}{z_1}+ \cdots + z_{k-i}\pfrac{}{z_{k-i-1}} \}, \mbox{ for } 0 \leq i \leq k-1,
$$
and
$$
\mathcal{C}^i =  \mspan \{\pfrac{}{z_{k-i+1}}, \cdots, \pfrac{}{z_k}\}, \; \dleq 1 i k-2.
$$
Notice that for each $0\leq i \leq k-3$, the only nontrivial condition for $[f+u_{0}^{i}g_{0}+u_{1}^{i}g_{1}, \mathcal{C}^{i+1}]\subset \mathcal{G}^{i}$ to be satisfied for $TCh_1^k$ is $[f+v_{0}^{i}g_{0}+v_{1}^{i}g_{1}, \frac{\partial}{\partial z_{k-i}}]\in \mathcal{G}^{i}$ implying $[f,\frac{\partial}{\partial z_{k-i}}]-v_{0}^{i} \frac{\partial}{\partial z_{k-i-1}} \in \mathcal{G}^{i}$ and hence
$$\frac{\partial f_{k-i-1}}{\partial z_{k-i}}(z)+v_{0}^{i}=0.$$
The latter is feedback invariant because $[f+u_{0}^{i}g_{0}+u_{1}^{i}g_{1}, \mathcal{C}^{i+1}]\subset \mathcal{G}^{i}$ is feedback invariant as explained just after the definition of $U_{sing}^{i}$ in Section~\ref{chap 3.1 sec:Flat-outputs}. Another argument proving feedback invariance is that we look for the vector field
$f(x)+u_{0}(x)^{i}g_{0}+u_{1}(x)^{i}g_{1}$ belonging to the affine distribution $f(x)+\mathcal{G}^{0}(x)$ which, obviously, is feedback invariant. To summarize, $v^{\ast}\in \bigcup_{i=0}^{k-3}U_{sing}^{i}(z^{\ast})$ if and only if
$$\frac{\partial f_{k-i-1}}{\partial z_{k-i}}(z^{\ast})+v_{0}^{\ast}=0, \, \dleq 0 i k-3.$$

To analyze the condition $[f+u_{0}^{k-2}g_{0}+u_{1}^{k-2}g_{1},l]\in \mathcal{G}^{k-2}$, where $l\in \mathcal{L}$ and $l \not\in \mathcal{C}^{k-2}$, take $l=\frac{\partial}{\partial z_{2}}$. Then
$$
[f+v_{0}^{k-2}g_{0}+v_{1}^{k-2}g_{1},l] = [f,\frac{\partial}{\partial z_{2}}]-v_{0}^{k-2}\frac{\partial}{\partial z_{1}}  \in \mathcal{G}^{k-2},
$$
if and only if 
$$
\frac{\partial f_{1}}{\partial z_{2}}(z)+v_{0}^{k-2}=0.
$$
The definition of $U_{_{\mathcal{L}-sing}}^{k-2}$ is feedback invariant (for the some reasons as those  giving invariance of $U_{sing}^{i}$, $0\leq i \leq k-3$) and thus $v^{\ast}\in U_{_{\mathcal{L}-sing}}^{k-2}$ if and only if $\frac{\partial f_{1}}{\partial z_{2}}(z^{\ast})+v_{0}^{\ast}=0$, {where $\mathcal{L}$ is such that $\mathcal{G}^0(x^*) \not\in \mathcal{L}(x^*)$. 
If $\mathcal{L}$ is such that $\mathcal{G}^0(x^*) \in \mathcal{L}(x^*)$, we will show when proving the equivalence $(i) \Longleftrightarrow (ii)$,  that under the assumption, which we always assume, $( d \varphi_{0} \wedge  d \varphi_{1} \wedge  d \dot{\varphi}_{0} \wedge d \dot{\varphi}_{1})(x^{\ast}, u^{\ast}) \neq~0$,  where $\mathcal{L}^{\bot}=\textrm{span} \{d\varphi_{0},d\varphi_{1}\}$, we have $u^{\ast}\not\in U_{_{\mathcal{L}-sing}}^{k-2}(x^*) $ and in $\mathcal{X}^*\times \mathbb{R}^2$, where $\mathcal{X}^* $ is a sufficiently small neighborhood of $x^*$, the set $U_{_{\mathcal{L}-sing}}^{k-2}(x) $ consists of two connected components that define, for each fixed value $x \in \mathcal{X}^*$, $x\neq x^*$, an affine subspace of $U= \mathbb{R}^2$.}

Now observe that the set of the singular control values $U_{_{\mathcal{L}-sing}}^{k-2}$ (at which {$(\vp_0, \vp_1)$ } ceases to be a flat output for $TCh_1^k$) is determined by $\mathcal{L}$ which, in turn, is uniquely associated to the choice of the flat output $(\varphi_{0},\varphi_{1})$ by $\mathcal{L}^{\bot}=\textrm{span} \{d\varphi_{0},d\varphi_{1}\}$. Different choices of $(\varphi_{0},\varphi_{1})$ lead, in general, to different distributions $\mathcal{L}$ and, consequently, to different singular control values and the system is not flat only at those that are singular for all choices of $\mathcal{L}$. Hence
$$
U_{sing}=\bigcup_{i=0}^{k-3}U_{sing}^{i} \cup U_{sing}^{k-2}
$$
where
$$
U_{sing}^{k-2}=\bigcap_{\mathcal{L}} U_{_{\mathcal{L}-sing}}^{k-2}.
$$

\medskip

\textit{Proof of (F2).}
{It was shown in \cite{li2010flat} that conditions $(FO2)$ and $(FO2)'$ are equivalent (for control-linear systems $\Sigma_{lin}$ but notice that $(FO2)$ and $(FO2)'$  do not involve the drift $f$). We deduce immediately that $(ii) \Leftrightarrow (iii) $. 
We will now prove that $(ii)\Rightarrow(i)$.}

First consider the case $(L_g\vp_0,L_g\vp_1)(x^*)\neq (0,0) $. By \cite{li2010flat}, a pair $(\varphi_{0},\varphi_{1})$ satisfying $(FO1)-(FO2)$ forms a flat output of the control-linear system $\Sigma_{lin}$ and, also by \cite{li2010flat}, $(\varphi_{0},\varphi_{1})$ is compatible with the chained form so there exists a local static feedback transformation bringing $\Sigma_{lin}$ into the chained form with $z_{0}=\varphi_{0}$ and $z_{1}=\varphi_{1}+k_{0}\varphi_{0}$, $k_{0}\in \mathbb{R}$, which thus transforms the control-affine system $\Sigma_{aff}$ into
$$
\begin{array}{l ccc c l}
 \dot z_0 =f_0(z)+v_0 & \dot z_1 &=& f_1(z) &+& z_2v_0\\
               &   & \vdots & & &\\
               &  \dot z_{k-1} &=& f_{k-1} (z) &+& z_kv_0\\
               &  \dot z_k &=& f_{k} (z) &+&v_1 
\end{array}
$$
Replacing $v_{0}$ by $v_{0}-f_{0}$ and $v_{1}$ by $v_{1}-f_{k}$ and using $[f,\mathcal{C}^{i}]\subset\mathcal{D}^{i}$, we conclude (repeating the proof of (F1)) that the system is in the triangular form and thus, flat at $(x^{\ast},u^{\ast})$ such that $u^{\ast}\not\in U_{_{\mathcal{L}-sing}}=\bigcup_{i=0}^{k-3}U_{sing}^{i} \cup U_{_{\mathcal{L}-sing}}^{k-2}$, where $\mathcal{L} = (\mspan\{d\vp_0,d\vp_1 \})^\perp$.
\medskip

Now consider the case $(L_g\varphi_0,\, L_g\varphi_1)(x^{\ast})= (0, 0)$.  Since $\Sigma_{aff}: \dot x = f(x) +u_0 g_0(x) +u_1 g_1(x) $ is locally, around $x^*$, feedback equivalent to $TCh_1^k$, we can assume that $\Sigma_{aff}$ is in the triangular form $TCh_1^k$ around $z^{\ast}=0$:
$$
 TCh_1^k\left\{
 \begin{array}{l lcl c l}
  \dot z_0 =v_0 & \dot z_1 &=&  f_1(z_0,z_1,z_2) &+& z_2v_0 \\
                & \dot z_2 &=&  f_2(z_0,z_1,z_2, z_3)&+& z_3v_0 \\
                &  &\vdots&  &&  \\
                &  \dot z_{k-1} &=& f_{k-1}(z_0,\cdots, z_k)&+&z_kv_0\\
                &   \dot z_{k} &=& v_1 & &
 \end{array}\right.
$$

The characteristic distribution $\mathcal{C}^{k-2}$ takes the form
$\mathcal{C}^{k-2} = \textrm{span}\, \{\frac{\displaystyle \partial}{\displaystyle \partial z_{3}},\dots,\frac{\displaystyle \partial}{\displaystyle \partial z_{k}}\},$
and the condition $ L_{c} \varphi_{i} =0$, for any $c\in \mathcal{C}^{k-2}$, given by {$(FO2)$ implies that $\varphi_{i}=\varphi_{i}(z_{0},z_{1},z_{2})$, for $i=0,1$.} 
Condition {$(FO1)$} implies that $ d \varphi_0 \wedge d \varphi_1 (x^{\ast}) \neq 0$, that is equivalent to
$$
\textrm{rk}
\left(
\begin{array}{c c c}
\displaystyle\frac{\partial \varphi_{0}}{\partial z_{0}} & \displaystyle\frac{\partial \varphi_{0}}{\partial z_{1}} & \displaystyle\frac{\partial \varphi_{0}}{\partial z_{2}}\\[2mm]
\displaystyle\frac{\partial \varphi_{1}}{\partial z_{0}} & \displaystyle\frac{\partial \varphi_{1}}{\partial z_{1}} & \displaystyle\frac{\partial \varphi_{1}}{\partial z_{2}}
\end{array}
\right)(0) =2.
$$
Notice that the condition $(L_g\varphi_0,\, L_g\varphi_1)(x^{\ast})= (0, 0)$ implies that $\frac{\partial \varphi_{0}}{\partial z_{0}}(0)=\frac{\partial \varphi_{1}}{\partial z_{0}}(0)=0$ and thus we get
$$
\textrm{rk}
\left(
\begin{array}{c c}
 \displaystyle\frac{\partial \varphi_{0}}{\partial z_{1}} & \displaystyle\frac{\partial \varphi_{0}}{\partial z_{2}}\\[2mm]
\displaystyle\frac{\partial \varphi_{1}}{\partial z_{1}} & \displaystyle\frac{\partial \varphi_{1}}{\partial z_{2}}
\end{array}
\right)(0) =2.
$$
We assume $\varphi_{0}(0) = \varphi_{1}(0)=0$ (if not, replace $\varphi_{0}$ by $\varphi_{0}-\varphi_{0}(0)$ and  $\varphi_{1}$ by $\varphi_{1}-\varphi_{1}(0)$).
We will introduce new coordinates $(\tilde{z}_{1},\tilde{z}_{2})=(\varphi_{0},\varphi_{1})$ in two steps. Assume that $\frac{\partial \varphi_{1}}{\partial z_{2}}(0) \neq 0$ (if not, permute $\varphi_{0}$ and $\varphi_{1}$) and put $\tilde{z}_{2}=\varphi_{1}(z_{0},z_{1},z_{2})$. Then the two first components become
\begin{eqnarray*}
 \dot z_1 &=&  \tilde{f}_1(z_0,z_1,\tilde{z}_2) + a(z_0,z_1,\tilde{z}_2)v_0 \\
 \dot{\tilde{z}}_{2}
 &=& \tilde{f}_{2}(z_{0},z_{1},\tilde{z}_{2},z_{3}) + b(z_{0},z_{1},\tilde{z}_{2},z_{3})v_{0},
\end{eqnarray*}
where $\tilde{f}_{2}=L_{f}\varphi_{1}$, $b=L_{g_0}\varphi_{1}$ and $a=z_{2}=\varphi_{1}^{-1}(z_{0},z_{1},\tilde{z}_{2})$ is the inverse of $\varphi_{1}$ with respect to~$z_{2}$. 
Notice that $b=L_{g_{0}}\varphi_{1}=\frac{\partial \varphi_{1}}{\partial z_{0}}+ \frac{\partial \varphi_{1}}{\partial z_{1}}z_{2}+ \frac{\partial \varphi_{1}}{\partial z_{2}}z_{3}$ is affine with respect to $z_{3}$ and $\frac{\partial \varphi_{1}}{\partial z_{2}}(0)\neq 0$ so $\tilde{z}_{i}=L_{g_{0}}^{i-3}b$, for $3\leq i \leq k$, is a valid local change of coordinates in which the system, under the feedback $\tilde{v}_{1}=L_{f}L_{g_{0}}^{k-3}b + v_{0}L_{g_{0}}^{k-2}b + v_{1}L_{g_{1}}L_{g_{0}}^{k-3}b$, takes the form
$$
 \begin{array}{l lcl c l}
  \dot z_0 =v_0 & \dot z_1 &=&  \tilde{f}_1(z_0,z_1,\tilde{z}_2) &+& a(z_0,z_1,\tilde{z}_2) v_0 \\
                & \dot{\tilde{z}}_2 &=&  \tilde{f}_2(z_0,z_1,\tilde{z}_2, \tilde{z}_3)&+& \tilde{z}_3v_0 \\
                &  &\vdots&  &&  \\
                &  \dot{\tilde{z}}_{k-1} &=& \tilde{f}_{k-1}(z_0,z_1,\tilde{z}_2,\cdots, \tilde{z}_k)&+&\tilde{z}_kv_0\\
                &   \dot{\tilde{z}}_{k} &=& \tilde{v}_1. & &
 \end{array}
$$
Now put $\tilde{z}_{1}=\varphi_{0}(z_0,z_1,z_{2})$. We get $\dot{\tilde{z}}_{1}=L_{f}\varphi_{0}+v_{0}L_{g_{0}}\varphi_{0}$. Notice that $L_{g_{0}}\varphi_{0}$ is affine with respect to $z_{3}$ and $L_{f}\varphi_{0}$ is, in general, nonlinear with respect to $z_{3}$ since so is $\tilde{f}_{2}$. Omitting $``\sim"$ we get
\begin{equation}\label{chap 3.1 Tri-form-A-B}
 \begin{array}{l lcl c l}
  \dot z_0 =v_0 & \dot{z}_1 &=&  f_1(z_0,z_1,z_2,z_3) &+& (A+Bz_{3}) v_0 \\
                & \dot{z}_2 &=&  f_2(z_0,z_1,z_2, z_3)&+&z_3v_0 \\
                &  &\vdots&  &&  \\
                &  \dot{z}_{k-1} &=& f_{k-1}(z_0,z_1,\cdots, z_k)&+&z_kv_0\\
                &   \dot{z}_{k} &=& v_1,& &
 \end{array}
\end{equation}
where $A$ and $B$ depend on $z_0,z_1,z_2$ only.
Observe that for (\ref{chap 3.1 Tri-form-A-B}), we have $\varphi_{0}=z_{1}$, $\varphi_{1}=z_{2}$ and $\mathcal{C}^{k-2} = \textrm{span}\, \{\frac{ \partial}{ \partial z_{3}},\dots,\frac{ \partial}{ \partial z_{k}}\}$, therefore the condition $(L_g\varphi_0) L_{[c,g]}\varphi_1 = (L_g\varphi_1) L_{[c,g]}\varphi_0$
gives $A+z_3B=z_3B$ and thus $A\equiv 0$ everywhere.

Notice that the function $f_2(z_0,z_1,z_2, z_3)$ can always be expressed as
$$ f_2(z_0,z_1,z_2, z_3)= f_{20}(z_0,z_1,z_2) + z_{3}f_{21}(z_0,z_1,z_2, z_3) $$
for some smooth functions $f_{20}$ and $f_{21}$ and thus
$$ \dot z_2 =  f_2(z_0,z_1,z_2, z_3)+ z_3v_0 =  f_{20}(z_0,z_1,z_2) + z_{3}(f_{21}(z_0,z_1,z_2, z_3)+v_{0}).$$
Define the new control $\tilde{v}_{0}=f_{21}(z_0,z_1,z_2, z_3)+v_{0}$ and denote $\eta = f_{21}$, then (\ref{chap 3.1 Tri-form-A-B}) becomes
\begin{equation}\label{chap 3.1 Tri-form-B}
 \begin{array}{l lcl c l}
  \dot z_0 =\tilde{v}_0 - \eta & \dot z_1 &=&  \tilde{f}_1(z_0,z_1,z_2,z_{3}) &+& z_{3}B\tilde{v}_0 \\
                & \dot z_2 &=&  \tilde{f}_2(z_0,z_1,z_2)&+& z_3\tilde{v}_0 \\
                &  &\vdots&  &&  \\
                &  \dot z_{k-1} &=& \tilde{f}_{k-1}(z_0,\cdots, z_k)&+&z_k\tilde{v}_0\\
                &   \dot z_{k} &=& v_1, & &
 \end{array}
\end{equation}
where $\tilde{f}_{2}=f_{20}$ and $\tilde{f}_{i}=f_{i}-z_{3}B\eta$, for $i\neq 2$.

Note that $\Sigma_{aff}$ is  assumed to be locally, around $x^* \in X$, static feedback equivalent to $TCh_1^k$, hence the conditions $[f,\mathcal{C}^{i}]\subset \mathcal{G}^{i}$ hold, for $1 \leq i \leq k-2$, and are invariant under change of coordinates and feedback. Clearly, for (\ref{chap 3.1 Tri-form-B}), $\mathcal{C}^{k-2} = \textrm{span}\, \{\frac{ \partial}{\partial z_{3}},\dots,\frac{ \partial}{ \partial z_{k}}\}$ and thus $[\tilde{f},\mathcal{C}^{k-2}]\subset \mathcal{G}^{k-2}$ implies $[\tilde{f},\frac{\partial}{\partial z_{3}}] \in \mathcal{G}^{k-2}$ and yields
$$
\left[\tilde{f},\frac{\partial}{\partial z_{3}}\right]=\left(
                   \begin{array}{c}
                     -\frac{\partial\eta}{\partial z_{3}}\\
                      \frac{\partial\tilde{f}_{1}}{\partial z_{3}}\\
                      0\\
                      \end{array}
                \right)
         = \alpha \left(
                   \begin{array}{c}
                     1\\
                      z_{3}B\\
                      z_{3}\\
                      \end{array}
                \right)+
                 \beta  \left(
                   \begin{array}{c}
                     0\\
                      B\\
                      1\\
                      \end{array}
                \right),
$$
modulo $\mathcal{C}^{k-2}$, for some smooth functions $\alpha,\beta$ which gives $ \frac{\partial\tilde{f}_{1}}{\partial z_{3}} =0$.
Therefore $\tilde{f}_1 = \tilde{f}_1(z_0,z_1,z_2)$ and thus (\ref{chap 3.1 Tri-form-B}) is, actually, in the following form
\begin{equation}\label{chap 3.1 Tri-form-B-3}
 \begin{array}{l lcl c l}
  \dot z_0 =\tilde{v}_0 - \eta \quad  & \dot z_1 &=&  \tilde{f}_1(z_0,z_1,z_2) &+& z_{3}B\tilde{v}_0 \\
                & \dot z_2 &=&  \tilde{f}_2(z_0,z_1,z_2)&+& z_3\tilde{v}_0 \\
                &  &\vdots&  &&  \\
                &  \dot z_{k-1} &=& \tilde{f}_{k-1}(z_0,\cdots, z_k)&+&z_k\tilde{v}_0\\
                &   \dot z_{k} &=& v_1, & &
 \end{array}
\end{equation}
with $(\varphi_{0},\varphi_{1})=(z_{1},z_{2})$. Define a new variable $y=z_{3}\tilde{v}_0$. Notice that, although $y=z_{3}\tilde{v}_{0}$ is not a valid control transformation (since $z_{3}^{\ast}=0$), it is a system's variable under the assumption that the differentials $dy=z_{3}d\tilde{v}_{0}+\tilde{v}_{0}dz_{3}$ is nonzero at $(z^{\ast},\tilde{v}_{0}^{\ast})$. Actually, $\dot{\varphi}_{0}$ and $\dot{\varphi}_{1}$ are functions of the system variables $z_{0},z_{1},z_{2}$ and $y$. Recall that $\varphi_{0}=z_{1}$ and $\varphi_{1}=z_{2}$.
The condition $ \textrm{rk} \frac{\partial (\varphi, \dot{\varphi})}{\partial(x,u)}(x^{\ast},u^{\ast})=4$ together with
$$
\frac{\partial(\varphi,\dot{\varphi})}{\partial(x,u)}= \frac{\partial(\varphi,\dot{\varphi})}{\partial(z_0,z_1,z_2,y)} \cdot \frac{\partial(z_0,z_1,z_2,y)}{\partial(x,u)}
$$
implies that $\textrm{rk} \frac{\partial (\dot{\varphi}_{0},\dot{\varphi}_{1})}{\partial (z_{0},y)}(z^{\ast},v^{\ast})=2$. By the implicit function theorem, we can express
\begin{eqnarray*}
z_{0}&=& \zeta_{0}(\varphi_{0},\varphi_{1},\dot{\varphi}_{0},\dot{\varphi}_{1})\\
y &=& \zeta_{y}(\varphi_{0},\varphi_{1},\dot{\varphi}_{0},\dot{\varphi}_{1})
\end{eqnarray*}
in a neighborhood of $(z^{\ast},v^{\ast})$, for some smooth functions $\zeta_{0},\zeta_{y}$.

We have $\dot z_0 = \t v_0 - \eta = v_0$ and $\dot z_2 = \t f_2 + z_3\t v_0 =\t f_2 + z_3(v_0+\eta)$. Recall that $\t f_2 $ depends on $z_0, z_1, z_2$ only. So knowing  $\dot z_0 =v_0$ and $\dot z_2 $, we can calculate $z_3$ using the implicit functions theorem if $v_0+\eta+ z_3 \pfrac{\eta}{z_3}\neq 0$. Then  {$\dot z_3$ gives $z_4$ if $v_0+\eta+ \pfrac{f_4}{z_4}\neq 0$ and so on, proving that indeed $(\vp_0,\vp_1)$ is an $x$-flat output at $(x^*, u^*)$. }

\medskip

{To conclude the proof, we have to show the implication $(i)\Rightarrow (ii)$.
When proving Proposition~\ref{chap 3.1 prop-charac-FO}, we will show that any flat output $(\varphi_{0},\varphi_{1})$ of a system $\Sigma_{aff}$ feedback equivalent to $TCh_1^k$ satisfies $(d\varphi_{0} \wedge d \varphi_{1} \wedge d\dot \varphi_{0} \wedge d \dot\varphi_{1})(x^{\ast},u^{\ast})\neq 0$ and
$L_c \varphi_0 =  L_c\varphi_1= (L_g \varphi_0) L_{[c,g]} \varphi_1 - (L_g \varphi_1) L_{[c,g]} \varphi_0 =  0$, for any $c \in \mathcal{C}^{k-2}$. 
If $(L_g\vp_0,L_g\vp_1)(x^*)\neq (0,0)$, we conclude in the same way as for item $(F1)$ that the singular control values $v^*$  coincide with $v^*\in U_{\mathcal{L}-sing}(z^*)$.} 

{Let us consider the case  $(L_g\vp_0,L_g\vp_1)(x^*)= (0,0)$. Since the conditions $L_c \varphi_0 =  L_c\varphi_1= (L_g \varphi_0) L_{[c,g]} \varphi_1 - (L_g \varphi_1) L_{[c,g]} \varphi_0 =  0$ are valid everywhere on $X$, we repeat the proof of  $(ii)\Rightarrow (i)$ and bring the system into the form (\ref{chap 3.1 Tri-form-B-3}), around $z^*=0$, with $(\vp_0, \vp_1)=(z_1,z_2).$}
Now we will show that the singular control values $v^*$ at which the procedures of calculating $z_0$ and $v_0$ fail, given by $\mrk \pfrac{(\dot \vp_0, \dot \vp_1)}{z_0,y}(z^*, v^*) \leq 1 $ and $v_0^* = -(\eta+z_3 \pfrac{\eta}{z_3})(z^*) $, coincide with $v^*\in U_{\mathcal{L}-sing}^{k-2}(z^*) $ and $v^*\in U_{sing}^{k-3}(z^*)  $, respectively.

To this end, calculate $U_{_{\mathcal{L}-sing}}^{k-2}(z)= \{v(z) = (v_{0},v_{1})^{\top}:[f+v_{0}g_{0}+v_{1}g_{1},l]\in \mathcal{G}^{k-2}\}$. 
Since $d\varphi_{0}=dz_{1}$ and $d\varphi_{1}=dz_{2}$, we have $\mathcal{L} = (\mspan \{d\varphi_0,d\varphi_1\})^{\bot} = \textrm{span} \, \{\frac{\partial}{\partial z_{0}}, \frac{\partial}{\partial z_{3}}, \frac{\partial}{\partial z_{4}}, \ldots,\frac{\partial}{\partial z_{k}}\}$ and $\mathcal{G}^{k-2}=\mathcal{L}+ \textrm{span}\,\{ B\frac{\partial}{\partial z_{1}}+\frac{\partial}{\partial z_{2}}\}$. 
Thus $[f+v_{0}g_{0}+v_{1}g_{1},l]\in \mathcal{G}^{k-2}$, for any $l\in \mathcal{L}$, holds (taking the only nontrivial case $l=\frac{\partial}{\partial z_{0}}$) if and only if 
$[f,\frac{\partial}{\partial z_{0}}]+v_{0}[g_{0},\frac{\partial}{\partial z_{0}}] \in \mathcal{G}^{k-2}$
which is equivalent to 
$[(\frac{\partial f_{1}}{\partial z_{0}} + v_{0}z_{3} \frac{\partial B}{\partial z_{0}})\frac{\partial}{\partial z_{1}} + \frac{\partial f_{2}}{\partial z_{0}}\frac{\partial}{\partial z_{2}}]\in \mathcal{G}^{k-2}$ and thus to $[(\frac{\partial f_{1}}{\partial z_{0}} + v_{0}z_{3} \frac{\partial B}{\partial z_{0}})\frac{\partial}{\partial z_{1}} + \frac{\partial f_{2}}{\partial z_{0}}\frac{\partial}{\partial z_{2}}]\wedge (B\frac{\partial}{\partial z_{1}}+\frac{\partial}{\partial z_{2}})=0$. 
This yields $v^{\ast} \in U_{\mathcal{L}-sing}^{k-2}(z^{\ast})$ if and only if $\frac{\partial f_{1}}{\partial z_{0}}(z^{\ast}) - B \frac{\partial f_{2}}{\partial z_{0}}(z^{\ast}) + v_{0}^{\ast}z_{3}^{\ast} \frac{\partial B}{\partial z_{0}}(z^{\ast}) =0$ which coincides with $\textrm{rk} \frac{\partial (\dot{\varphi}_{0},\dot{\varphi}_{1})}{\partial (z_{0},y)}(z^{\ast},v^{\ast})\leq 1$.

{Notice that under the assumption $( d \varphi_{0} \wedge  d \varphi_{1} \wedge  d \dot{\varphi}_{0} \wedge d \dot{\varphi}_{1})(z^{\ast}, u^{\ast}) \neq 0$, we have $\frac{\partial f_{1}}{\partial z_{0}}(z^{\ast}) - B \frac{\partial f_{2}}{\partial z_{0}}(z^{\ast})\neq 0$ and, since $z^{\ast} =0 $, it follows that $v_{0}^{\ast}\not\in U_{_{\mathcal{L}-sing}}^{k-2}(z^*) $. Moreover, since $\frac{\partial B}{\partial z_{0} }\neq 0$ (otherwise $\mathcal{G}^{k-1}\neq TX$), for each fixed value $x\neq x^*$ in $\mathcal{X}^*$, a sufficiently small neighborhood of $x^*$, we get $(v_0,v_1) \in U_{_{\mathcal{L}-sing}}^{k-2}(z^*)$ with $v_0 = \frac{\psi(z_0,z_1,z_2)}{z_3 }$, where $\psi =  (\frac{\partial f_{1}}{\partial z_{0}} )(\frac{\partial B}{\partial z_{0} })^{-1}$, and $v_1$ any. Thus  in $\mathcal{X}^*\times \mathbb{R}^2$,  the set $U_{_{\mathcal{L}-sing}}^{k-2}(x) $ consists of two connected components that define, for each fixed value $x \in \mathcal{X}^*$, $x\neq x^*$, an affine 
subspace of $U= \mathbb{R}^2$.}

To analyze $v_0^* = -(\eta+z_3 \pfrac{\eta}{z_3})(z^*) $, notice that for (\ref{chap 3.1 Tri-form-B-3}), $\mathcal{C}^{k-2} = \mspan\{\pfrac{}{z_3}, \cdots, \pfrac{}{z_n}\}$ and $\mathcal{G}^{n-3} = \mathcal{C}^{k-2} +\mspan\{\pfrac{}{z_0} + z_3 B \pfrac{}{z_1} +z_3\pfrac{}{z_2} \}$. It follows that $[\t f+\t v_0 \t g_0+ \t v_1\t g_1, \mathcal{C}^{k-2} ] \in \mathcal{G}^{n-3}$ is equivalent to $[\t f+\t v_0 g_0+ \t v_1 g_1, \pfrac{}{z_3}] \wedge (\pfrac{}{z_0} + z_3 B \pfrac{}{z_1} +z_3\pfrac{}{z_2} ) = 0 \mmod \mathcal{C}^{k-2}$, which yields $-\pfrac{\eta}{z_3} + \t v_0 (\pfrac{}{z_1} +z_3\pfrac{}{z_2})\wedge(\pfrac{}{z_0} + z_3 (B \pfrac{}{z_1} +z_3\pfrac{}{z_2})) ) = 0$ implying $z_3 \pfrac{\eta}{z_3} + \t v_0 =  z_3 \pfrac{\eta}{z_3} + \eta+ v_0= 0$. Thus, indeed, $v_0^*= -(z_3 \pfrac{\eta}{z_3} + \eta)(z^*)$ if and only if  $v^*\in U_{sing}^{n-3}(z^*)$.
\end{proof}

{\subsection{Proof of Proposition \ref{chap 3.1 prop output m=1 S_lin}}}
\begin{proof}
 In \cite{li2010flat}, the equivalence of the following conditions has been proven for {any two-input system feedback equivalent to the chained form and for a pair of smooth functions $(\vp_0, \vp_1)$:}
\begin{enumerate}[${(}i{)}$]
\item The pair $(\vp_0,\vp_1)$ is an $x$-flat output of $\S_{lin}$ at $(x^*,u^*)$, where $u^*$ is such that $u^*_{0}g_0(x^*)+u^*_{1}g_1(x^*)\not\in\mathcal{C}^{1}(x^*)$;
\item The pair $(\vp_0,\vp_1)$  satisfies the following conditions:
    \begin{enumerate}[${(FO}1{_{lin})}$]
		\item $d\varphi_0 \wedge d\varphi_1(x^*) \neq 0$;
		\item  $L_c \varphi_0 =  L_c\varphi_1=  L_c (\frac{L_g\varphi_1}{L_g\varphi_0})=  0,$ for any $c \in \mathcal C^{k-2}$, where the functions $\varphi_0, \,\varphi_1$ are ordered such that $L_g\varphi_0(x^*) \neq 0$, which is always possible due to item $(FO3_{lin})$;
		\item $(L_g\varphi_0(x^*),\, L_g\varphi_1(x^*))\neq (0,\, 0)$;
    \end{enumerate}
 \item The pair $(\vp_0,\vp_1)$  satisfies the following conditions:
    \begin{enumerate}[${(FO}1{_{lin})'}$]
		\item $d\varphi_0 \wedge d\varphi_1(x^*) \neq 0$;
		\item $\mathcal{L} = (\mspan \{d\varphi_0,d\varphi_1\})^\perp\subset \mathcal{G}^{k-2}$;
		\item $\mathcal{G}^0(x^*) \not\subset\mathcal{L}(x^*)$.
    \end{enumerate}   
\end{enumerate}

{In the view of the above, item $(F3)$ is obvious. So is  $(F6)$ because $(FO1)'$ yields  $(FO1_{lin})'$, the condition $(L_g\varphi_0(x^*),\, L_g\varphi_1(x^*))\neq (0,\, 0)$ implies $(FO3_{lin})'$, and  $(FO2)'$ and $(FO2_{lin})'$ coincide.}

{
To show $(F5)$, notice that $(FO2)'$ and $(FO2_{lin})'$ coincide. 
To prove that $(\vp_0, \vp_1)$ satisfies $(F01)$, we can bring, see \cite{li2010flat}, the control-linear system $\S_{lin}$ into the chained form compatible with the flat output $(\vp_0, \vp_1)$ (which is assumed to be a flat output of  $\S_{lin}$), that is, $Ch_1^k$ with $z_0 =\vp_0$ and $z_1= \vp_1$. 
In the $z$-coordinates, the drift takes the triangular form for $TCh_1^k$. 
By a direct calculation, we can check that $(d\vp_0 \wedge d\vp_1 \wedge d\dot\vp_0 \wedge d\dot\vp_1)(z^*,v^*)\neq 0,$ where $v^*\not\in U_{\mathcal{L}-sing}(z^*)$ and $\mathcal{L} = (\mspan\{d\vp_0, d\vp_1 \})^\perp$.
Hence  $(\vp_0, \vp_1)$ is an $x$-flat output of  $\S_{aff}$ at $(x^*,\tilde u^*)$ where $\tilde u^*\not\in  U_{\mathcal{L}-sing}(x^*)$.
}

{It remains to prove $(F4)$. If $(\vp_0, \vp_1)$ is a flat output of  $\S_{lin}$, then the conditions $(FO1_{lin})-(FO3_{lin})$
are satisfied and thus so are $(FO1)-(FO2)$ because $(FO2)$ and $(FO2_{lin})$ coincide and $(\vp_0, \vp_1)$ being a flat output of $\S_{lin}$ satisfies $(FO1)$ with $\dot \vp_i = L_{F_{lin}} \vp_i$, $i=0,1$.}

{To prove the converse}, we have to show that condition $(F01)$ $( d \varphi_{0} \wedge  d \varphi_{1} \wedge  d \dot{\varphi}_{0} \wedge d \dot{\varphi}_{1})(x^{\ast}, u^*) \neq~0$, where  $\dot{\varphi}_{i}$, for $i=0,1$ is understood as $\dot{\varphi}_{i}=L_{F_{lin}}\varphi_{i}$ and $F_{lin}=u_0g_0+u_1g_1$, implies that
$(L_g\vp_0,L_g\vp_1)(x^*)\neq (0,0)$. 

Bring $\S_{lin}$ into the chained form $Ch_1^k $ around $z^* =0$ and let $(\vp_0, \vp_1)$ be a flat output. Since $L_c \varphi_0 =  L_c\varphi_1=0$, for all $c \in \mathcal C^{k-2} = \mspan\{\pfrac{}{z_3}, \cdots,\pfrac{}{z_k} \}$, it follows $\varphi_i = \varphi_i(z_0,z_1,z_2) $, for $i=0,1$. 
Assume $(L_g\vp_0,L_g\vp_1)(0)=(0,0)$, otherwise the claim holds. 
Thus  $\pfrac{\vp_i}{z_0}(0) = 0$, for $i=0,1$, and since $( d \varphi_{0} \wedge  d \varphi_{1} )(0) \neq 0$, we deduce $\mrk \pfrac{(\vp_0, \vp_1)}{(z_1,z_2)}(0) =2$. 
Assume that $\frac{\partial \varphi_{1}}{\partial z_{2}}(0) \neq 0$ (if not, permute~$\varphi_{0}$ and $\varphi_{1}$) and put $\tilde{z}_{2}=\varphi_{1}$. 
Notice that $b=L_{g_{0}}\varphi_{1}=\frac{\partial \varphi_{1}}{\partial z_{0}}+ \frac{\partial \varphi_{1}}{\partial z_{1}}z_{2}+ \frac{\partial \varphi_{1}}{\partial z_{2}}z_{3}$ is affine with respect to~$z_{3}$ and $\frac{\partial \varphi_{1}}{\partial z_{2}}(0)\neq 0$ so $\tilde{z}_{i}=L_{g_{0}}^{i-3}b$, for $3\leq i \leq k$, is a valid local change of coordinates in which the system, under the feedback $\tilde{v}_{1}=v_{0}L_{g_{0}}^{k-2}b + v_{1}L_{g_{1}}L_{g_{0}}^{k-3}b$, takes the form
$$
 \begin{array}{l lcl c l}
  \dot z_0 =v_0 & \dot z_1 &=& a(z_0,z_1,\tilde{z}_2) v_0 \\
                & \dot{\tilde{z}}_2 &=&  \tilde{z}_3v_0 \\
                &  &\vdots&  &&  \\
                &  \dot{\tilde{z}}_{k-1} &=& \tilde{z}_kv_0\\
                &   \dot{\tilde{z}}_{k} &=& \tilde{v}_1. & &
 \end{array}
$$
where  $ a = z_2 = \vp_1^{-1}(z_0,z_1,\t z_2)$. 
The condition $(L_g\varphi_0)L_{[c,g]}\varphi_1= (L_g\varphi_1)L_{[c,g]}\varphi_0$ yields $\pfrac{\vp_0}{z_0} + a \pfrac{\vp_0}{z_1} =~0$.
So omitting the tildes, we obtain $\dot \vp_0 = \pfrac{\vp_0}{z_2} {z}_3v_0 = \pfrac{\vp_0}{z_2} \dot \vp_1$. 
Therefore  the differentials satisfy $d\dot \vp_0 = \dot \vp_1 d\pfrac{\vp_0}{z_2} \mmod \mspan\{d\dot \vp_1\}$ and since $\dot \vp_1(0) = 0 $, we get $(d\dot \vp_0\wedge d\dot \vp_1)(0) =0$, which contradicts the independence of flat outputs and their differentials. Thus $(L_g\vp_0,L_g\vp_1)(0)\neq(0,0)$.
{Now it is obvious that $L_c (\frac{L_g\varphi_1}{L_g\varphi_0})=  0$ is equivalent to $(L_g\varphi_0)L_{[c,g]}\varphi_1= (L_g\varphi_1)L_{[c,g]}\varphi_0$, where $L_g\varphi_0(x^*) \neq 0$ (after permuting $ \vp_0$ and $ \vp_1$, if necessary)}. 
\end{proof}
{\subsection{Proof of Proposition \ref{chap 3.1 prop:output-eq:m=1}}}
\begin{proof}
For the proof of Proposition~\ref{chap 3.1 prop:output-eq:m=1} in the case  $L_{g}\vp_0(x^*) \neq 0$, we refer the reader to~\cite{li2010flat}. Let us consider the case $L_{g}\vp_0(x^*) = 0$. Bring the system $\S_{aff}$ into the form $TCh_1^k$, around $z^*=0$. The characteristic distribution $\mathcal{C}^{k-2}$ takes the form
$\mathcal{C}^{k-2} = \textrm{span}\, \{\frac{\displaystyle \partial}{\displaystyle \partial z_{3}},\dots,\frac{\displaystyle \partial}{\displaystyle \partial z_{k}}\},$
and the condition $ L_{c} \varphi_{0} =0$, for any $c\in \mathcal{C}^{k-2}$, implies that $\varphi_{0}=\varphi_{0}(z_{0},z_{1},z_{2})$.
From $<d\vp_0, \mathcal{G}^{k-2}>(0)\neq 0$, we deduce $\pfrac{\vp_0}{z_2}(0)\neq 0$. Introducing the new coordinate $\t z_2 = \vp_0$ and following exactly the proof of item $(F2)$ of Theorem~\ref{chap 3.1 thm:output:m=1}, we get (omitting the tildes for $\t z$)
\begin{equation}
 \begin{array}{l lcl c l}
  \dot z_0 =\tilde{v}_0 - \eta(z_0,z_1,z_2,z_3) & \dot z_1 &=&\t {f}_1(z_0,z_1,z_2,z_3) &+&a(z_0,z_1,z_2)\tilde{v}_0 \\
                & \dot z_2 &=&  \tilde{f}_2(z_0,z_1,z_2)&+& z_3\tilde{v}_0 \\
                &  &\vdots&  &&  \\
                &  \dot z_{k-1} &=& \tilde{f}_{k-1}(z_0,\cdots, z_k)&+&z_k\tilde{v}_0\\
                &   \dot z_{k} &=& v_1, & &
 \end{array}
\end{equation}
with $\varphi_{0}=z_{2}$.
The condition $[f, \mathcal{C}^{k-2}]\in  \mathcal{G}^{k-2} $ implies $\pfrac{f_1}{z_3} = -a \pfrac{\eta}{z_3} $.  In these coordinates we have $v = (L_{g}\vp_0)[c_{k-2},g]-(L_{[c_{k-2},g]}\vp_0)g$ $= z_3\pfrac{}{z_2} -(\pfrac{}{z_0}+a\pfrac{}{z_1}+z_3\pfrac{}{z_2})\mmod \mathcal{C}^{k-2}.$ {The} distribution $\mathcal{L} =  \mathcal{C}^{k-2}+\mspan\{\pfrac{}{z_0}+a\pfrac{}{z_1}\}$ is, {indeed}, involutive and of corank two in $TX$. {Thus} there exists 
 a smooth function $ \psi = \psi(z_0,z_1,z_2)$ such that $\pfrac{\psi}{z_1}(0)\neq 0$ and $\pfrac{\psi}{z_0}+a \pfrac{\psi}{z_1} = 0$ and we put $\t z_1 =\psi $. Then $\dot {\t z}_1 = L_f\psi + \pfrac{\psi}{z_2}z_3\tilde{v}_0 =$ $\bar {f}_1(z_0,z_1,z_2,z_3) + z_3B(z_0,z_1,z_2)\tilde{v}_0$. From  $[f, \mathcal{C}^{k-2}]\in  \mathcal{G}^{k-2} $, it follows that $\bar {f}_1 = \bar {f}_1(z_0,z_1,z_2)$. We have 
$$
 \begin{array}{l lcl c l}
  \dot z_0 =\tilde{v}_0 - \eta & \dot {\t z}_1 &=&  \bar{f}_1(z_0,z_1,z_2) &+& z_{3}B\tilde{v}_0 \\
                & \dot z_2 &=&  \tilde{f}_2(z_0,z_1,z_2)&+& z_3\tilde{v}_0 \\
                &  &\vdots&  &&  \\
                &  \dot z_{k-1} &=& \tilde{f}_{k-1}(z_0,\cdots, z_k)&+&z_k\tilde{v}_0\\
                &   \dot z_{k} &=& v_1, & &
 \end{array}
$$
with $\psi = \t z_1$ and $\varphi_{0}=z_{2}$. The pair $(\vp_0, \psi)=$ $ (z_2, z_1) $ is an $x$-flat output at $(z^*, v^*)$, with  $v^*\not\in U_{\mathcal{L}-sing}(z^*)$, if and only if $(\frac{\partial \bar f_{1}}{\partial z_{0}}- B \frac{\partial \t f_{2}}{\partial z_{0}})(0)\neq 0$, i.e., $(d \psi \wedge d \dot \psi \wedge d \varphi_{0} \wedge d \dot \varphi_{0})(0)\neq 0$. 
\end{proof}

\subsection{Proof of Proposition \ref{chap 3.1 prop-charac-FO}}

\begin{proof}
Consider $\Sigma_{aff}$ static feedback equivalent to $TCh_1^k$ and let $(\varphi_{0},\varphi_{1})$ be a flat output {at $(x^{\ast}, u^{\ast})$, such that $(L_{g}\vp_0,L_{g}\vp_1)(x^*) \neq (0,0)$,  where $g$ is an arbitrary vector field in $\mathcal{G}$ such that $g(x^*)\not\in \mathcal{C}^{k-2}(x^*)$}. 
Form the decoupling matrix $D=(D_{ij})$, where $D_{ij}=L_{g_{j}}\varphi_{i}$, $0\leq i,j\leq1$. The involutive closure $\bar{\mathcal{G}}^{0}$ of $\mathcal{G}^{0}$ is $TX$, so $1\leq \textrm{rk}\, D(x) \leq 2$.
If $\textrm{rk}\, D(x)=2$, then via a suitable feedback transformation $\dot{\varphi}_{i}=\tilde{v}$, $i=0,1$, which contradicts flatness. 
Thus $\mrk D(x)=1$ {in a neighborhood of $x^*$}, since $(L_{g}\vp_0,L_{g}\vp_1)(x^*) \neq (0,0)$.
We have $d\varphi_{0} \wedge d\varphi_{1}(x)\neq 0$ so put $z_{0}=\varphi_{0}$, $z_{1}=\varphi_{1}$ and, after applying feedback, the first two components of the transformed system $\dot{z}=f+v_{0}g_{0}+v_{1}g_{1}$ become $\dot{z}_{0}=v_{0}$, $\dot{z}_{1}=a_{1}(z)+b_{1}(z)v_{0}$. 
The successive time-derivatives $\varphi_{1}^{(l)}$ of $\varphi_{1}=z_{1}$ cannot depend on $v_{1}$, for $0 \leq l \leq k-1$ (it would contradict flatness) and the $k$-th derivative depends explicitly on $v_{1}$, {otherwise 
we would obtain a contradiction with the independence of flat outputs and their time-derivatives at $(x^*,u^*)$}. 
Notice, however, that $\varphi_{1}^{(l)}$ is a polynomial of degree $l$, with respect to~$v_{0}$, with the leading coefficient being $L_{g_{0}}^{l-1}b_{1}$. 
Since $\varphi_{1}^{(l)}$ does not depend on $v_1$, for $ 1 \leq l \leq k-1$, it follows that $
L_{g_1}L_{g_0}^{l-1}b_1=0$  for $\dleq 1 l k-2$. 
We claim that the functions  $z_0,$ $ z_1,$ $ b_{1},$ $\ldots, $ $L_{g_{0}}^{k-2}b_{1}$ are independent {at any point of an open and dense $X'\subset X$. If not, take $x_0$ and its open neighborhood $V \subset X \backslash X'$ and} let 
$s$ be the {largest} integer such that $z_0,$ $ z_1,$ $ b_{1},$ $\ldots, $ $L_{g_{0}}^{s}b_{1}$ are independent  {in $V$}.
Assume $s\leq k-3$. Introduce new coordinates $z_i = L_{g_0}^{i-2}b_1$ {in $V$}, for $\dleq 2 i s$. We get: 
$$
\begin{array}{llcl@{\ }c@{\ }l l}
\dot z_0 =v_0 & \dot z_1 &=& a_1(z) &+& z_{2}v_0\\
	      & \dot z_2 &=& a_2(z) &+& z_{3}v_0\\
               &   & \vdots & & &\\
               &  \dot z_{s+1} &=& a_{s+1} (z) &+& z_{s+2}v_0\\
               &  \dot z_{s+2} &=& a_{s+2} (z) &+& b_{s+2} (z_{0},\ldots,z_{s+2})v_0\\               
               &  \dot{\bar z} &=&\bar f &+& \bar g_0 v_0& + \bar g_1 v_1 
\end{array}
$$
where $\bar z = (z_{s+3}, \ldots,z_k)$. 
Notice that the vector field $[g_0, g_1]$ is of the form $\sum_{i=s+3}^k \a_i\pfrac{}{z_i}$, with~$\a_i$ smooth functions. We deduce that $\bar{ \mathcal{G}}^0$, the involutive closure of $\mathcal{G}^0 = \mspan\{g_0, g_1\}$, satisfies  $\bar{ \mathcal{G}}^0 \subset \mspan\{g_0, \pfrac{}{z_{s+3}}, \cdots, \pfrac{}{z_{k}}\}$. This yields $\bar{ \mathcal{G}}^0\neq TX$, which contradicts the fact that for $\S_{aff},$ static feedback equivalent to $TCh_1^k$, we have $\bar{ \mathcal{G}}^0=TX$. Thus $s = k-2$ and we put
%
$
z_{2}=b_{1},\ldots,z_{k}=L_{g_{0}}^{k-2}b_{1},
$
and replace $v_{1}$ by
{$ L_{f}L_{g_{0}}^{k-2}b_{1}+v_{0}(L_{g_{0}}^{k-1}b_{1}) + v_{1}(L_{g_{1}}L_{g_{0}}^{k-2}b_{1})$}. We get
$$
g_{0}=\frac{\partial}{\partial z_{0}}+z_{1}\frac{\partial}{\partial z_{2}}+\cdots+z_{k-1}\frac{\partial}{\partial z_{k}} \quad \textrm{and} \quad g_{1}=\frac{\partial}{\partial z_{k}}.
$$
Using exactly the same arguments as in sufficiency part of the proof of Theorem~\ref{chap 3.1 thm:m=1} (the forms of $\mathcal{G}^{i}$ and of $\mathcal{C}^{i}$ and the condition $[f,\mathcal{C}^{i}]\in \mathcal{G}^{i})$ we conclude that on {$X^{\prime}$}, open and dense in $X$, the system is {locally} in the triangular form
\begin{equation*}
TCh_1^k:
\left\{
\begin{array}{l lcl c l}
 \dot z_0 =v_0 & \dot z_1 &=& f_1(z_{0},z_{1},z_{2}) &+& z_{2}v_0\\
               &   & \vdots & & &\\
               &  \dot z_{k-1} &=& f_{k-1} (z_{0},\ldots,z_{k}) &+& z_{k}v_0\\
               &  \dot z_k &=& v_1 & &
\end{array}
\right.
\end{equation*}
The flat output $(\varphi_{0},\varphi_{1})=(z_{0},z_{1})$ satisfies
$$
L_c \varphi_0 =  L_c\varphi_1= (L_g \varphi_0) L_{[c,g]} \varphi_1 - (L_g \varphi_1) L_{[c,g]} \varphi_0 =0,
$$
where $c \in \mathcal{C}^{k-2}=\textrm{span}\,\{\frac{\partial}{\partial z_{3}},\ldots,\frac{\partial}{\partial z_{k}}\}$ and $g$ is any vector field such that $\mathcal{G}^{0}=\textrm{span}\,\{g,c_{1}\}$ where $c_{1}=\frac{\partial}{\partial z_{k}}$ is the characteristic vector field of $\mathcal{G}^{1}$. 
In order to prove that we can bring the system into the triangular form {$TCh_1^k$,}  around any $x^{\ast}\in X$ (and not only on~{$X^{\prime}$}), notice that the characteristic distribution $\mathcal{C}^{k-2}$ is defined everywhere (not only on~{$X^{\prime}$})  so, by continuity, the conditions $L_c \varphi_0 =  L_c\varphi_1=(L_g \varphi_0) L_{[c,g]} \varphi_1 - (L_g \varphi_1) L_{[c,g]} \varphi_0=0$ hold everywhere on $X$ implying that if we put the control system $\Sigma_{aff}$, around an arbitrary point $x^{\ast}\in X$, into the triangular form {$TCh_1^k$,}  then for the flat output $(\varphi_{0},\varphi_{1})$, we have
$\varphi_{i}=\varphi_{i}(z_{0},z_{1},z_{2})$, $0\leq i \leq 1$, {on $X'$ and thus on $X$}.

{Since we have assumed that $(L_{g}\vp_0,L_{g}\vp_1)(x^*) \neq (0,0)$, we can apply the following change of coordinates (permute $\vp_0$ and $\vp_1$, if necessary)
$z_0 = \vp_0$, $z_1= \vp_1$ and $z_i ={L_{g_0}^{i-2}}\psi$, for $\dleq 2 i k$, where $\psi = \frac{L_{g_0}\vp_1}{L_{g_0}\vp_0} $, in which the control vector fields are in the chained form  with $(\vp_0, \vp_1) = (z_0, z_1)$.
The system $\S_{aff}$ is assumed to be feedback equivalent to the triangular form $TCh_1^k$, hence satisfies the compatibility condition \textit{(Comp)}. 
Using the $z$-coordinates and applying the feedback $f \mapsto f-(L_f\vp_0) g_0 - (L_f^{k-1}\psi) g_1$, we transform $\S_{aff}$  into the triangular form $TCh_1^k$ with $(\varphi_{0},\varphi_{1})=(\t z_{0},\t z_{1})$ around any $x^{\ast}\in X$.}

Notice that we have proved, in particular, that any flat output $(\varphi_{0},\varphi_{1})$ of a system $\Sigma_{aff}$ feedback equivalent to $TCh_1^k$ satisfies $(d\varphi_{0} \wedge d \varphi_{1} \wedge d\dot \varphi_{0} \wedge d \dot\varphi_{1})(x^{\ast},u^{\ast})\neq 0$ and $L_c \varphi_0 =  L_c\varphi_1= (L_g \varphi_0) L_{[c,g]} \varphi_1 - (L_g \varphi_1) L_{[c,g]} \varphi_0 =  0$, for any $c \in \mathcal{C}^{k-2}$, that is, conditions $(FO1)-(FO2)$ of Theorem~\ref{chap 3.1 thm:output:m=1}.

\end{proof}

\subsection{Proof of Theorem \ref{chap 3.1 thm:flat_outputs:m>=2}}
\begin{proof}[Proof of (m-F1)]
Consider a control-affine system $\Sigma: \dot x = f(x) +\sum_{i=0}^m u_i g_i(x)$ locally, around~$x^*$, static feedback equivalent to $TCh_m^k$, and bring it into the form $TCh_m^k$, around~$z^*$. For simplicity of notation, we continue to denote by $f$, respectively by $g_i$, for $\dleq 0 i m$, the drift, respectively the controlled vector fields of  $TCh_m^k$.

It is clear that $TCh_m^k$ is $x$-flat, with $\varphi = (z_0,z_1^1, \cdots,z_m^1 )$ being a flat output, at any point  $(z^*,v^*) \in X\times \mathbb{R}^{m+1} $ satisfying
$$
\mrk F^l(z^*) = m, \mbox{ for } 1 \leq l \leq k-1,
$$
where $F^l$, for $ 1 \leq l \leq k-1,$ is the $m\times m$ matrix given by
$$
F^l_{ij} =\pfrac{(f^l_j + z^{l+1}_jv^*_{0})}{z^{l+1}_i} ,  \mbox{ for } 1 \leq i, j \leq m.
$$
Moreover, the differential weight of $\varphi = (z_0,z_1^1, \cdots,z_m^1 )$ is $(k+1)(m+1)$, since expressing $z$ and~$v$ involves $\vp_i^{(j)}$, for $\dleq 1 i m$ and $\dleq 0 j k$.

Recall that in coordinates $z$, using the notation $\mspan \{\pfrac{}{z^i}\} = \mspan \{\pfrac{}{z^i_1}, \cdots, \pfrac{}{z^i_m}\}$, we have
$$
\mathcal{G}^i = \mspan \{\pfrac{}{z^{k-i}}, \cdots, \pfrac{}{z^k}, g_0 \}, \; \dleq 0 i k-1,
$$
$$
\mathcal{C}^i =  \mspan \{\pfrac{}{z^{k-i+1}}, \cdots, \pfrac{}{z^k}\}, \; \dleq 1 i k-2,
$$
 and
 $$
 \mathcal{L}=  \mspan \{\pfrac{}{z^{2}}, \cdots, \pfrac{}{z^k}\}.
 $$
We have $\mathcal{C}^1=  \mspan \{\pfrac{}{z^{k}_1}, \cdots, \pfrac{}{z^k_1}\},$ and thus
$$
\begin{array}{lcl}
 \mathcal{G}^0+ [f+gv,\mathcal{C}^{1} ] &=& {\mathcal{G}^0+  \mspan \{ [f+gv,\pfrac{}{z^{k}_j}], \, \dleq 1 j m\}} \vspace{0.2cm}\\
 &=& \mathcal{G}^0+  \mspan \{ \pfrac{(f^{k-1}_1+z^k_1v_0)}{z^{k}_j}\pfrac{}{z^{k-1}_1}+ \cdots +\pfrac{(f^{k-1}_m+z^k_mv_0 )}{z^{k}_j }\pfrac{}{z^{k-1}_m}, \dleq 1 j m \},
  \end{array} 
$$
where $gv= \sum_{i=0}^m g_iv_i$.
By induction, we obtain
\begin{small}
$$
\mathcal{G}^i+ [f+gv,\mathcal{C}^{i+1} ]= $$
$$\mathcal{G}^i+ \mspan \{
\pfrac{(f^{k-i-1}_1+z^{k-i}_1v_0)}{z^{k-i}_j}\pfrac{}{z^{k-i-1}_1}+ \cdots +\pfrac{(f^{k-i-1}_m+z^{k-i}_mv_0)}{z^{k-i}_j}\pfrac{}{z^{k-i-1}_m}, \ \dleq 1 j m \}.
$$
\end{small}
{Therefore for any $\dleq 0 i k-2$, we have $\mrk F^{i+1}(z^*,v^*)=m$ if and only if $\mrk(\mathcal{G}^i+ [f+gv,\mathcal{C}^{i+1} ])(z^*,v^*) = (i+2)m+1$, for $\dleq 0 i k-3$, and $\mrk(\mathcal{G}^{k-2}+ [f+gv,\mathcal{L})(z^*,v^*) = km+1$. It follows that the original system $\S_{aff}$ is $x$-flat at $(x^*,u^*)$ such that $ u^* \not\in U_{m-sing}(x^*)$, of differential weight at most $(k+1)(m+1)$.} 

%
\bigskip

As we have noticed, $(\vp_0, \ldots,\vp_m) = (z_0, z^1_1,\ldots,z^1_m)$ is an $x$-flat output of $TCh_{m}^k$ of differential weight $(k+1)(m+1)$ since expressing $z$ and $v$ involves $\varphi_{i}^{(j)}$, for $0\leq j \leq k$.

Now, we will show (which is interesting as an independent observation) that the differential weight of any $x$-flat output of $\Sigma_{aff}: \dot{x}=f+\sum_{i=0}^{m}u_{i}g_{i}$, with $m+1$ controls and $km+1$ states, is at least $(k+1)(m+1)$. Let $(\vp_0, \ldots,\vp_{m})$ be an $x$-flat output of $\Sigma_{aff}$. Define $D=(D_{ij})$, where $D_{ij}=L_{g_{i}}\varphi_{j}$ and put $r(x)= \textrm{rk}\,D(x)$. Clearly, $r(x)$ is constant on an open and dense subset $X^{\prime}$ of $X$ (so denote it $r(x)=r$) and choose $x_{0}\in X^{\prime}$. By a suitable (local) change of coordinates and static invertible feedback, we get
\begin{equation*}
\begin{array}{l ccc c l}
{ \dot z^0} =v^0 \qquad& { \dot z^1} &=& A^{1}(z) &+& B^{1}(z)v^0\\
               &  { \dot z^2} &=& A^{2}(z) &+& B^{2}(z)v
\end{array}
\end{equation*}
where $\textrm{dim}\, z^{0}=r$, $\textrm{dim}\, z^{1}=m-r+1$, $z_{0}^{0}=\varphi_{0},\ldots,z_{r-1}^{0}= \varphi_{r-1}$ and $z_{r}^{1}=\varphi_{r} ,\ldots,z_{m}^{1}=\varphi_{m} $.


Due to flatness we can express (with the help of the flat outputs $\varphi_{i}$ and their time-derivatives) $mk+1$ components of $z$ and $m+1$ components of $v$, i.e., $m(k+1)+2$ functions. 
Using $\varphi_{i}=z_{i}^{0}$ and $\dot{\varphi}_{i}=v_{i}^{0}$, $0\leq i \leq r-1$, we express $2r$ system variables. 
The remaining $m(k+1)+2-2r$ system variables (that is,  the components of $z^{1}$, $z^{2}$ and the remaining components of $v$) depend on derivatives of $\varphi_{i}$, $r\leq i \leq m$.
Denote by $s_{i}$ the maximal order of the derivative $\varphi_{i}^{(s_{i})}$, $r \leq i \leq m$, that is involved. Put $s=\textrm{max}\{s_{i}: r\leq i \leq m\}$. By taking the time-derivatives of $\varphi_{i}$ up to order  $s_{i}\leq s$, we can express at most $(s+1)(m-r+1)$ functions. This number cannot thus be smaller than the number of functions that remain to be expressed, that is, we need
$$
(s+1)(m-r+1) \geq m(k+1)+2-2r,
$$
which is equivalent to
$$
m(s-k) \geq (r-1)(s-1).
$$
Now, three cases are possible. It is clear that if $s<k$, then the left hand side is negative, so the inequality is not satisfied. If $s=k$, then either $r=1$ or $s=1$. The latter is impossible since $s\geq 2$.
In the case $r=1$, we have $\textrm{dim}\, z^{0}= \textrm{dim}\, v^{0} =1$ and 
in order to express all $m(k+1)+2$ variables of the system, we will use $s=k$ derivatives $v^{0}$, $\dot{v}^{0},\ddot{v}^{0},\ldots,(v^{0})^{(s-1)}$. Thus the differential weight of~$\varphi$  is at least $m(k+1)+s+1= m(k+1)+k+1 = (m+1)(k+1)$.

Finally, if $s>k$, then there exists $\varphi_{j}$, for some $r+1 \leq j \leq m+1$, that we differentiate $s$ times so it involves at least $s-1$ time derivatives of $\dot{\varphi_{j}}=A_{j}^{1}(z) +B_{j}^{1}(z)v^{0}$, where $A_{j}^{1}$ is the $j$-th component of $A^{1}$ and $B_{j}^{1}$ is the $j$-th row of $B^{1}$. 
The involutive closure $\bar{\mathcal{G}}^{0}$ of the distribution $\mathcal{G}^{0}$ is $TX$ so $B_{j}^{1}$ is nonzero.
It implies that $\varphi_{j}^{(s)}$ depends nontrivially on (at least one) component of $(v^{0})^{(s-1)}$. To summarize, we use $mk+1$ functions to express $z$, $m+1$ functions to express $v$, and we also use the $s-1$ derivatives $\dot{v}^{0},\ddot{v}^{0},\ldots,(v^{0})^{(s-1)}$, which gives at least $(k+1)(m+1)+1$ functions (since $s > k$). Therefore the differential weight is higher than $(k+1)(m+1)$ on $X^{\prime}$ and thus on $X$.

It remains to prove that the differential weight of any flat output (not necessary an $x$-flat output) cannot be smaller than $(k+1)(m+1)$. Let $(\vp_0, \ldots, \vp_m)$ be  an $(x,$ $u,$ $\dot{u},$ $\ldots,u^{(p)})$-flat output of $\Sigma_{aff}$. 
Denote by $s_i$ the highest derivative of $\vp_i$, for $\dleq 0 i m$, involved in expressing the state $x$ and the control $u$, that is, by flatness, $\mathcal{ X} + \mathcal{U} \subset \Phi$, where $\mathcal{ X} = \mspan\{dx_1, \cdots, dx_n\}$, $\mathcal{ U} = \mspan\{du_0, \cdots, du_m\}$ and $\Phi = \mspan\{d\vp_i^{(j_i)}, \dleq 0 i m, \dleq{0}{j_i}{s_i}\}$. 
Let $s_{i^*}$ be the largest among the integers $s_i$.  
Either~$\vp_{i^*}$ depends on $u^{(l)}$, with $l\geq 1$ (but not on derivatives of $u$ higher than $l$) or $\vp_{i^*}$ depends on~$u$ (but not on derivatives of $u$) or $\vp_{i^*}$ depends on $x$ only. Then the differentials $\vp_{i^*}^{(j)}$ are independent modulo $\mathcal{ X} + \mathcal{U}$, for $\dleq{0}{j}{s_{i^*}}$ (in the first case),  for $\dleq{1}{j}{s_{i^*}}$ (in the second case) and for 
$\dleq{2}{j}{s_{i^*}}$ (in the third case, since $\dot\vp_{i^*}$ depends 
on $u$ because $\bar{\mathcal{G}}^0 = TX$). It follows that  $\mathcal{ X} + \mathcal{U} \subset \Psi =\mspan\{ d\vp_{i^*}, d\dot\vp_{i^*}, d\vp_i^{(j_i)}, \dleq 0 {i} m,i\neq i^*, \dleq{0}{j_i}{s_i}\}$.

We claim that $s_{i^*}\geq k $. If not, then $s_i\leq s_{i^*}\leq k-1$, for $\dleq 0 i m$ (recall that $s_{i^*} =\mathrm{max}\{s_i : \dleq 0 i m\}$),which implies $\mrk \Psi \leq mk+2 < m(k+1)+2 = \mrk(\mathcal{ X} + \mathcal{U})$, contradicting  $\mathcal{ X} + \mathcal{U} \subset \Psi$. Thus $s_{i^*}\geq k $. 

We have $\mathcal{ X} + \mathcal{U} \subset \Phi$ (by flatness) and $d\ddot\vp_{i^*}, \cdots,$ $ d\vp_{i^*}^{(s_{i^*})}$ belong to $\Phi$ and are independent modulo $\mathcal{ X} + \mathcal{U}$, so $\mrk \Phi \geq \mrk(\mathcal{ X} + \mathcal{U})+k-1 = $ $m(k+1)+2+k-1=$ $(m+1)(k+1)$ proving that the differential weight of $\vp$ is at least $(m+1)(k+1)$. 
Notice that $\mrk \Phi = (m+1)(k+1)$ if and only if $ s_{i^*} = s_i =k$, for any $\dleq 0 i m$, implying that with $\vp_i$, $i\neq i^*$, we express $mk$ system variables and the remaining two variables are expressed with  $\vp_{i^*}$. We deduce immediately that, in this case, all $\vp_i$ depend on $x$ only. 

\bigskip

\textit{Proof of (m-F2).} Let $(\vp_0, \cdots,\vp_m)$ be a minimal $x$-flat output for $\S_{aff}$. When proving \textit{ (m-F1)} we have shown that we can bring the system into the form 
\begin{equation*}
\begin{array}{l ccc c l}
{ \dot z_0} =v_0 & { \dot z^1} &=& A^{1}(z) &+& B^{1}(z)v^0\\
               &  { \dot z^2} &=& A^{2}(z) &+& B^{2}(z)v
\end{array}
\end{equation*}
where $z_{0} =\varphi_{0}$ and $z_{1}^{1}=\varphi_{1},\ldots,z_{m}^{1}=\varphi_m$ and  $\textrm{dim}\, z_0 =\textrm{dim}\, v_0 = 1$, being a consequence of the minimal differential weight $(k+1)(m+1)$ of $\vp$. 
For $\dleq i i m$, denote by $k_i$ the minimal integer such that $\vp_i^{(k_i)}$ depends explicitly on at least one $v_j$, for $\dleq 1 j m$. Since $\S_{aff}$ is static feedback equivalent to $TCh_m^k$, it follows that $k_i\leq k$. In order to prove that $k_i=k$, for $\dleq 1 j m$, suppose that there exists $k_i< k$ and assume, for simplicity, that $k_1< k$. 
Denote $\vp_1^{(k_1)} = v_1$ (with $v_1$ depending on $v_0, \cdots v_0^{(k_1-1)})$. 

Like in the the proof of \textit{ (m-F1)}, notice that 
due to flatness we can express (with the help of the flat outputs $\varphi_{i}$ and their time-derivatives) $mk+1$ components of $z$ and $m+1$ components of $v$, i.e., $m(k+1)+2$ functions. 
Using $\varphi_{0}=z_{0}$ and  $\varphi_{1}=z_{1}^{1}$,   we can express $2 + k_1+1 = k_1+3$  variables of the system. 
The remaining $m(k+1)+2-(k_1+3)$ system variables depend on derivatives of~$\varphi_{i}$, $2\leq i \leq m$.
Denote by $s_{i}$ the maximal order of the derivative $\varphi_{i}^{(s_{i})}$, $2\leq i \leq m$, that is involved. Put $s=\textrm{max}\{s_{i}: 2\leq i \leq m\}$. By taking the time-derivatives of $\varphi_{i}$ up to order  $s_{i}\leq s$, we can express at most $(s+1)(m-1)$ functions. This number cannot thus be smaller than the number of functions that remain to be expressed, that is, we need
$$
(s+1)(m-1)\geq m(k+1)+2-(k_1+3),
$$
which is equivalent to
$$
m(s-k) \geq s-k_1.
$$
We have $k_1<k$ so the inequality can be satisfied only if $s>k$, but this give the differential weight of $\vp$ at least $m(k+1)+2+s-1 \geq (k+1)(m+1)+2$, implying that $\vp$ is not a minimal flat output. It follows that for all $\dleq 1 i m$ we must have $ k_i =k$ (and the inequality is satisfied only in this case). The distribution $\mathcal{L}$ $=(\mspan\{d\vp_0, \cdots,d\vp_m \})^\perp$ is involutive (as annihilator of exact 1-forms) and satisfies $\mathcal{L}\subset \mathcal{G}^{k-2}$ (because all $ k_i =k$), as well as  $ \mathcal{G}^{0}(x^*)\not\subset\mathcal{L}(x^*)$ (since $g_0(x^*)\not\in \mathcal{L}(x^*)$).
It follows that $\mathcal{G}^0$ is in the $m$-chained form in $z$-coordinates, where $z_0 = \vp_0$, $z_i^j = L_{g_0}^{j-1} \vp_i$, for $\dleq 1 i m,$  $\dleq 1 j k$ (see Appendix B). The compatibility condition \textit{(m-Comp)} implies that  $\S_{aff}$ is in the triangular form.

\bigskip

\textit{Proof of (m-F3).} We will prove the implications: $(i) \Rightarrow (iii) \Rightarrow (ii) \Rightarrow (i)$.

$(i) \Rightarrow (iii)$.
Assume that the system $\S_{aff}: \dot x = f(x) +\sum_{i=0}^m u_ig_i(x)$ is $x$-flat at $(x^*,u^*)$,  where $u^* \not\in U_{m-sing}(x^*)$, and let $(\vp_0, \cdots,\vp_m)$ be its minimal $x$-flat output defined in a neighborhood $\mathcal{X}^*$ of $x^*$. It is well known that the differentials of flat outputs are independent at $x^*$, {thus implying \textit{(m-FO1)}.
By item  \textit{(m-F2)}, that we have just proven,} we can bring $\S_{aff}$, around any point {$x \in {\mathcal{X}}^*$} into the triangular form compatible with the chained form $TCh_m^k$,  with $(\vp_0, \cdots,\vp_m) = (z_0, z^1_1,\cdots,z^1_m)$  and $x^*$ transformed into $z^*\in \mathbb{R}^{km+1}$. 
In coordinates $z$, the corank one involutive subdistribution $\mathcal{L}$ of $\mathcal{G}^{k-2}$ is given by
$$
\mathcal{L}=  \mspan \{\pfrac{}{z^{2}}, \cdots, \pfrac{}{z^k}\},
$$
because it is unique and we immediately have
$$
\mathcal{L}^\perp=  \mspan \{d\vp_0, \cdots,d\vp_m\},
$$
which gives \textit{(m-FO2)} on ${\mathcal{X}^*}$. 
\bigskip

$(iii) \Rightarrow (ii)$. Suppose that the $(m+1)$-tuple $(\vp_0, \cdots,\vp_m)$  fulfills conditions {\textit{(m-FO1)-(m-FO2)}}. We apply the change of coordinates and the invertible feedback transformation presented in Appendix B 
(with $\phi_i$ replaced by $\vp_i$ and $\t u$ by $v$) that bring the control-linear system $\S_{lin}: \dot x = \sum_{i=0}^m u_ig_i(x)$ into the $m$-chained form, with $z_0 = \vp_0$ and $z_i^1 = \vp_i$, for $\dleq{1}{i}{m}.$ Thus $(\vp_0, \cdots,\vp_m) = (z_0, z_1^1, \cdots, z_m^1)$ is a  minimal $x$-flat output of $Ch_m^k$ at any $(z^*,v^*)$, with $v^*\neq 0$. It follows that $(\vp_0, \cdots,\vp_m)$ is a minimal $x$-flat output of $\S_{lin}$ at any $(x^*,\tilde u^*)$, with $\tilde u^*$ such that $\sum_{i=0}^m \tilde  u^*_{i}g_i(x^*)\not\in\mathcal{C}^{1}(x^*)$.

\bigskip

$(ii) \Rightarrow (i)$. Assume that the system $\S_{lin}: \dot x = \sum_{i=0}^m u_ig_i(x)$ is $x$-flat at $(x^*,\tilde u^*)$,  where $\tilde u^*$ is such that  $\sum_{i=0}^m \tilde u^*_{i}g_i(x^*)\not\in\mathcal{C}^{1}(x^*)$, where $\mathcal{C}^{1}$ is the characteristic distribution of $\mathcal{G}^1$. 
Let $(\vp_0, \cdots,\vp_m)$ be its minimal $x$-flat output defined in a neighborhood $\mathcal{X}$ of $x^*$. 
It is known, see \cite{li2011geometry}, that the minimal flat output satisfies $\mathcal{L}^\perp = \mspan\{d\vp_0, \cdots,d\vp_m\}$. 
By the construction given in Appendix B, bring  the system into the  $m$-chained form $Ch_m^k$ such that $(\vp_0, \cdots,\vp_m) = (z_0, z_1^1, \cdots,z_m^1)$ and $z^j_i ={L_{g_0}^{j-2}}\psi_i$, for $\dleq 2 j k$ and $\dleq 1 i m$, where $\psi_i = \frac{L_{g_0}\vp_i}{L_{g_0}\vp_0} $.
The system $\S_{aff}$ is assumed to be feedback equivalent to the triangular form $TCh_m^k$, hence satisfies the compatibility condition \textit{(m-Comp)}. 
Using the $z$-coordinates and applying the feedback $f \mapsto f-\sum_{i=0}^{m}\a_ig_i$, where $\a_0 = L_f\vp_0$ and $\a_i = L_f^{k-1}\psi_i$, we transform $\S_{aff}$   into the triangular form $TCh_m^k$.
We have proved, when showing \textit{(m-F1)}, that $(\vp_0, \cdots,\vp_m) = (z_0, z_1^1, \cdots,z_m^1)$ is an $x$-flat output of  $\S_{aff}$  at $(x^*,u^*)$ such that $ u^*\not\in U_{m-sing}(x^*).$
 \end{proof}

\section*{Appendices}
 \addcontentsline{toc}{section}{Appendices}

\subsection*{A. Involutive subdistribution of corank one}  
Consider  a non involutive distribution $\mathcal{G}$ of rank $d$, defined on a manifold $X$ of dimension $n$
and define its annihilator $\mathcal{G}^{\perp}=\{\omega\in \Lambda^1(X)\,:\,<\omega,f >=0, \forall f\in\mathcal{G}\}$.
Let $\omega_1, \dots, \omega_s$, where $s= n-d$, be differential 1-forms locally spanning the annihilator of $\mathcal{G}$, that is $\mathcal{G}^{\perp}=\mathcal{I} = \mspan \{\omega_1, \dots, \omega_s\}$.  The \emph{Engel rank} of $\mathcal{G}$ equals 1 at $x$ if and only if $(d\omega_i\wedge d\omega_j)(x) =0 \mmod \mathcal{I} ,$ for any $1\leq i,j\leq s$. For any $\omega \in \mathcal{I}$, we define $\mathcal{W}(\omega) = \{f \in \mathcal{G}: f \lrcorner d\omega\in \mathcal{G}^\perp\}$, where $\lrcorner$ is the interior product. The characteristic distribution $\mathcal{C} = \{f \in \mathcal{G}: [f,\mathcal{G}]\subset \mathcal{G}\}$ of $\mathcal{G}$ is given by
$$
\mathcal{C} = {\bigcap}_{i=1}^s \mathcal{W}(\omega_i).
$$
It follows directly from the Jacobi identity that the characteristic distribution is always involutive. Let $\mrk [\mathcal{G},\mathcal{G}]= d+r$.
Choose the differential forms $\omega_1, \dots, \omega_r, \dots,\omega_s$ such that $\mathcal{I} = \mspan \{\omega_1, \dots, \omega_s\}$ and $\mathcal{I}^1 = \mspan \{\omega_{r+1}, \dots, \omega_s\}$, where $\mathcal{I}^1$ is the annihilator of $[\mathcal G, \mathcal G]$. Define the distribution
$$
\cB = \sum_{i=1}^r\mathcal{W}(\omega_i).
$$
We have the following result proved by \cite{bryant1979some}, see also \cite{pasillas2001contact}.

\begin{proposition}
Consider a distribution $\mathcal G$ of rank $d$ and let $\mrk [\mathcal{G},\mathcal{G}]= d+r$.
\begin{enumerate}[(i)]
 \item Assume $r\geq 3$. The distribution $\mathcal G$ contains an involutive subdistribution of corank one if and only if it satisfies
\begin{enumerate}[{(ISD}1)]
 \item The Engel rank of $\mathcal G$ equals one;

 \item The characteristic distribution $\mathcal{C}$ of $\mathcal{G}$ has rank $d-r-1$.
\end{enumerate}
\noindent
Moreover, that involutive subdistribution is unique and is given by $\cB$.

 \item Assume $r= 2$. The distribution $\mathcal G$ contains a corank one subdistribution $\mathcal L$ satisfying $[\mathcal L, \mathcal L]\subset~\mathcal G$ if and only it verifies (ISD1)-(ISD2). In that case, $\cB$ is the unique distribution with the desired properties.

 \item Assume $r= 1$.  The distribution $\mathcal G$ contains an involutive subdistribution of corank one if and only it satisfies the condition (ISD2). In the case $r=1$, if an involutive subdistribution of corank one exists, it is never unique.
\end{enumerate}
\label{chap 3.1 prop_Bryant}
\end{proposition}

\subsection*{B. Constructing coordinates for the $m$-chained form } 
In \cite{pasillas2001canonical}, the following characterization of the  $m$-chained form was stated and proved:
 An $(m+1)$-input driftless control system $\S_{lin}: \dot x = \sum_{i=0}^mu_ig_i(x)$, with $m\geq 2$, defined on a manifold $X$ of dimension $km+1$, is locally static feedback equivalent, in a small neighborhood of a point $x^*\in X$,  to the $m$-chained form if and only if its associated distribution $\mathcal{G} = \mspan\{g_0, \cdots, g_m\}$ satisfies conditions \textit{(m-Ch1)-(m-Ch3)} of Theorem~\ref{chap 3.1 thm:m>=2}.

 The prove of this result provides a method to compute the diffeomorphism bringing any control system, for which it is possible, to the $m$-chained form. Now, we will explain how to do it.

 The involutive subdistribution $\mathcal{L}$ is unique and can be explicitly calculated (see Appendix A). 
 Choose $m+1$ independent functions $\phi_0$, $\phi^1_1, \cdots,$ $\phi^1_m$ whose differentials annihilates $\mathcal{L}$, that is
 $$\mspan\{d\phi_0, d\phi^1_1, \cdots, d\phi^1_m\} = (\mathcal{L})^\perp,$$
 and a vector field $g\in \mathcal{G}^0$ (which always exists due to condition $(m-Ch3)$) such that  $g(x^*)\not\in \mathcal{L}^{k-2}(x^*)$. Without loss of generality, we can assume $g = g_0$ and $L_{g_0}\phi^0_0(x^*)\neq 0$ (otherwise permute the vector fields $g_i$ or the functions $\phi^1_i$). Define the coordinates
 $$
 \left\{
\begin{array}{lcl}
 z_0&= &\phi_0\\
 z^1_i&= &\phi^1_i, \; \dleq 1 i m,\\
 z^j_i& = &\phi^j_i = \frac{L_{g_0}\phi^{j-1}_i}{L_{g_0}\phi_0},\; \dleq 1 i m, \; \dleq 2 j k,
\end{array}
\right.
$$
and the feedback
 $$
 \t u_0  =  u_0L_{g_0}\phi_0 \, \mbox{ and }\, \t u_j  = \sum_{i=0}^m u_iL_{g_i}\phi_j^k, \; \dleq 1 j m. $$
In the above coordinates, the distribution $\mathcal{G}$ takes the form
$$
\mathcal{G} = \mspan\{\pfrac{}{z^k_1}, \cdots,\pfrac{}{z^k_m}, \pfrac{}{z_0}+\sum_{j=1}^m \sum_{i=1}^{k-1}z^{i+1}_j\pfrac{}{z_j^i} \}
$$
and, equivalently, $\S_{lin}$ takes the $m$-chained form.


\bibliographystyle{apalike}
\bibliography{TRI-IJC}
\end{document}